\documentclass{article}
\usepackage{amsmath, epsfig,amssymb, multicol}
\usepackage[active]{srcltx}
\newtheorem{lemma}{Lemma}
\newtheorem{theorem}{Theorem}

\newcommand{\C}{{\cal C}}

\newcommand{\G}{{\cal G}}

\title{A Fractional Analogue of Brooks' Theorem}
\author{
Andrew D. King
\thanks{Simon Fraser University, Burnaby, British Columbia, Canada
({\tt andrew.d.king@gmail.com}).  This author was supported in part by an EBCO/Ebbich Postdoctoral Scholarship and the NSERC Discovery Grants of Pavol Hell and Bojan Mohar.}
\and
Linyuan Lu
\thanks{University of South Carolina, Columbia, SC 29208,
({\tt lu@math.sc.edu}). This author was supported in part by NSF
grant  DMS 1000475. }
  \and Xing Peng
\thanks{University of South Carolina, Columbia, SC 29208,
({\tt pengx@mailbox.sc.edu}). This author was supported in part by
NSF grant  DMS 1000475. } }
\begin{document}
\maketitle
\begin{abstract}
  Let $\Delta(G)$ be the maximum degree of a graph $G$.
  Brooks' theorem states that the only connected graphs with chromatic
  number $\chi(G)=\Delta(G)+1$ are complete graphs and odd cycles.  We
  prove a fractional analogue of Brooks' theorem in this
  paper.  Namely, we classify all  connected graphs $G$ such that the
  fractional chromatic number $\chi_f(G)$ is at least $\Delta(G)$.
  These graphs are complete graphs, odd cycles, $C^2_8$, $C_5\boxtimes K_2$,
  and graphs whose clique number $\omega(G)$ equals the maximum degree
  $\Delta(G)$.  Among the two sporadic graphs, the graph $C^2_8$ is
  the square graph of cycle $C_8$ while the other graph $C_5\boxtimes K_2$
  is the strong product of $C_5$ and $K_2$.  In fact, we prove a
  stronger result; if a connected graph $G$ with $\Delta(G)\geq 4$ is
  not one of the graphs listed above, then we have $\chi_f(G)\leq
  \Delta(G)- \frac{2}{67}$.
\end{abstract}
\section{Introduction}
The chromatic number of graphs with bounded degrees has been studied
for many years.  Brooks' theorem perhaps is one of the most
fundamental results; it is included by many textbooks on
graph theory. Given a simple connected graph $G$, let $\Delta(G)$ be
the maximum degree, $\omega(G)$ be the clique number,
 and $\chi(G)$ be the chromatic number.  Brooks'
theorem states that $\chi(G)\leq \Delta(G)$ unless $G$ is a complete
graph or an odd cycle.
 Reed \cite{reed1} proved that $\chi(G)\leq
\Delta(G)-1$ if $\omega(G)\leq \Delta(G)-1$ and $\Delta(G)\geq \Delta_0$
for some large constant $\Delta_0$. This excellent result was proved by
probabilistic methods, and $\Delta_0$ is at least
hundreds. Before this result, Borodin  and Kostochka \cite{bk} made
the following conjecture.

{\bf Conjecture \cite{bk}:}
Suppose that $G$ is a connected graph.
If $\omega(G)\leq \Delta(G)-1$
and $\Delta(G)\geq 9$, then we have
$$\chi(G)\leq \Delta(G)-1.$$

\begin{figure}[htbp]
\label{fig:1}
 \centerline{\psfig{figure=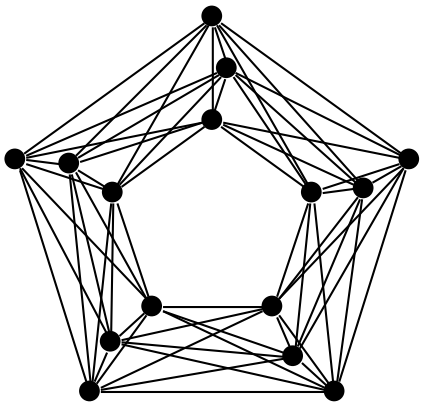, width=0.4\textwidth} }
\caption{The graph $C_5 \boxtimes K_3$.}
\end{figure}

If the conjecture is true, then it is best possible since there is a
$K_8$-free graph $G=C_5\boxtimes K_3$ (actually $K_7$-free, see Figure \ref{fig:1}) with
$\Delta(G)=8$ and $\chi(G)=8$.

Here we use the following notation of the strong product. Given two graphs $G$
and $H$, the {\it strong product} $G\boxtimes H$ is the graph with vertex
set $V(G)\times V(H)$,  and  $(a,x)$ is connected to $(b,y)$ if one of the following holds
\begin{itemize}
\item $a=b$ and $xy\in E(H)$,
\item $ab\in E(G)$ and $x=y$,
\item $ab\in E(G)$ and $xy\in E(H)$.
\end{itemize}

 Reed's result \cite{reed1}  settled Borodin and Kostochka's conjecture for sufficiently large $\Delta(G)$,
but the cases with small $\Delta(G)$ are hard to cover using the probabilistic method.

In this paper we  consider a fractional analogue of this problem.
The fractional chromatic number $\chi_f(G)$ can be defined as
follows. A $b$-{\it fold coloring} of $G$ assigns a set of $b$
colors to each vertex such that any two adjacent vertices receive
disjoint sets of colors. We say a graph $G$ is $a$:$b$-{\it
colorable} if there is a $b$-fold coloring of $G$ in which each
color is drawn from a palette of $a$ colors. We refer to such a
coloring as an $a$:$b$-coloring. The $b$-{\it fold coloring number},
denoted by $\chi_b(G)$, is the smallest integer $a$ such that $G$
has an $a$:$b$-coloring. Note that $\chi_1(G)=\chi(G)$.  It was shown
that $\chi_{a+b}(G) \leq \chi_a(G) + \chi_b(G) $.  The {\it
  fractional chromatic number} $\chi_f(G)$ is
 $ \underset{b \rightarrow \infty} \lim \frac{\chi_b(G)}{b}.$

By the definition, we have $\chi_f(G)\leq \chi(G)$.
 The fractional chromatic number can be viewed as a relaxation of the
 chromatic number. Many problems involving the chromatic number can be asked
again using the  fractional chromatic number. The fractional analogue often
has a simpler solution than the original
problem. For example, the famous $\omega-\Delta-\chi$ conjecture of
 Reed \cite{reed} states that for any simple graph $G$, we have
$$\chi(G)\leq \left\lceil \frac{\omega(G)+\Delta(G)+1}{2}\right
\rceil.$$
The fractional analogue of $\omega-\Delta-\chi$ conjecture was
proved by Molloy and Reed \cite{mr}; they actually proved a stronger
result with ceiling removed, i.e.,
\begin{equation} \label{eq:1}
\chi_f(G)\leq  \frac{\omega(G)+\Delta(G)+1}{2}.
\end{equation}

In this paper, we classify all connected graphs $G$ with
$\chi_f(G)\geq \Delta(G)$.

\begin{theorem}\label{main}
  A connected graph $G$ satisfies $\chi_f(G)\geq \Delta(G)$ if and only if
$G$ is one of the following
\begin{enumerate}
\item a complete graph,
\item  an odd cycle,
\item  a graph with  $\omega(G)=\Delta(G)$,
\item  $C^2_8$,
\item  $C_5\boxtimes K_2$.
\end{enumerate}
\end{theorem}

For the complete graph $K_n$, we have $\chi_f(K_n)=n$ and
$\Delta(K_n)=n-1$.  For the odd cycle $C_{2k+1}$, we have
$\chi_f(C_{2k+1})=2+\frac{1}{k}$ and $\Delta(C_{2k+1})=2$. If $G$ is
neither a complete graph nor an odd cycle but contains a clique of
size $\Delta(G)$, then we have
\begin{equation} \label{less}
\Delta(G)\leq \omega(G)\leq \chi_f(G)\leq \chi(G)\leq \Delta(G).
\end{equation}
The last inequality is from Brooks' theorem. The sequence of
inequalities above implies $\chi_f(G)=\Delta(G)$.

If $G$ is a vertex-transitive graph,
then we have \cite{su}
 $$\chi_f(G)=\frac{|V(G)|}{\alpha(G)},$$
where $\alpha(G)$ is the independence number of $G$.
Note  that both graphs $C^2_8$ and $C_5\boxtimes K_2$
are vertex-transitive and have  the independence number $2$.
Thus we have
$$\chi_f(C^2_8)=4=\Delta(C^2_8)\quad \mbox { and }\quad
\chi_f(C_5\boxtimes K_2)=5= \Delta(C_5\boxtimes K_2).$$

\begin{figure}[htbp]
 \centerline{ \psfig{figure=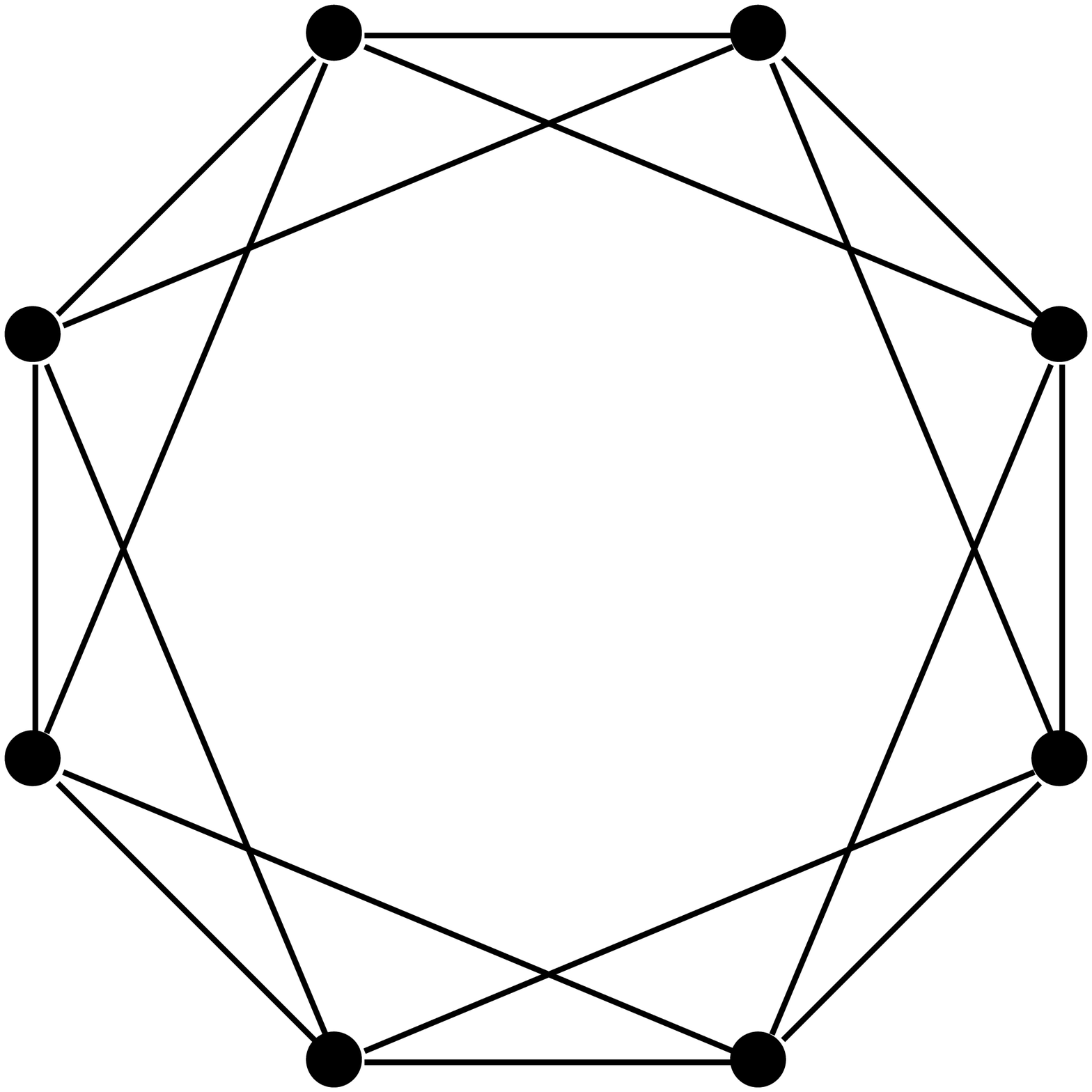, width=0.24\textwidth} \hspace{4cm} \psfig{figure=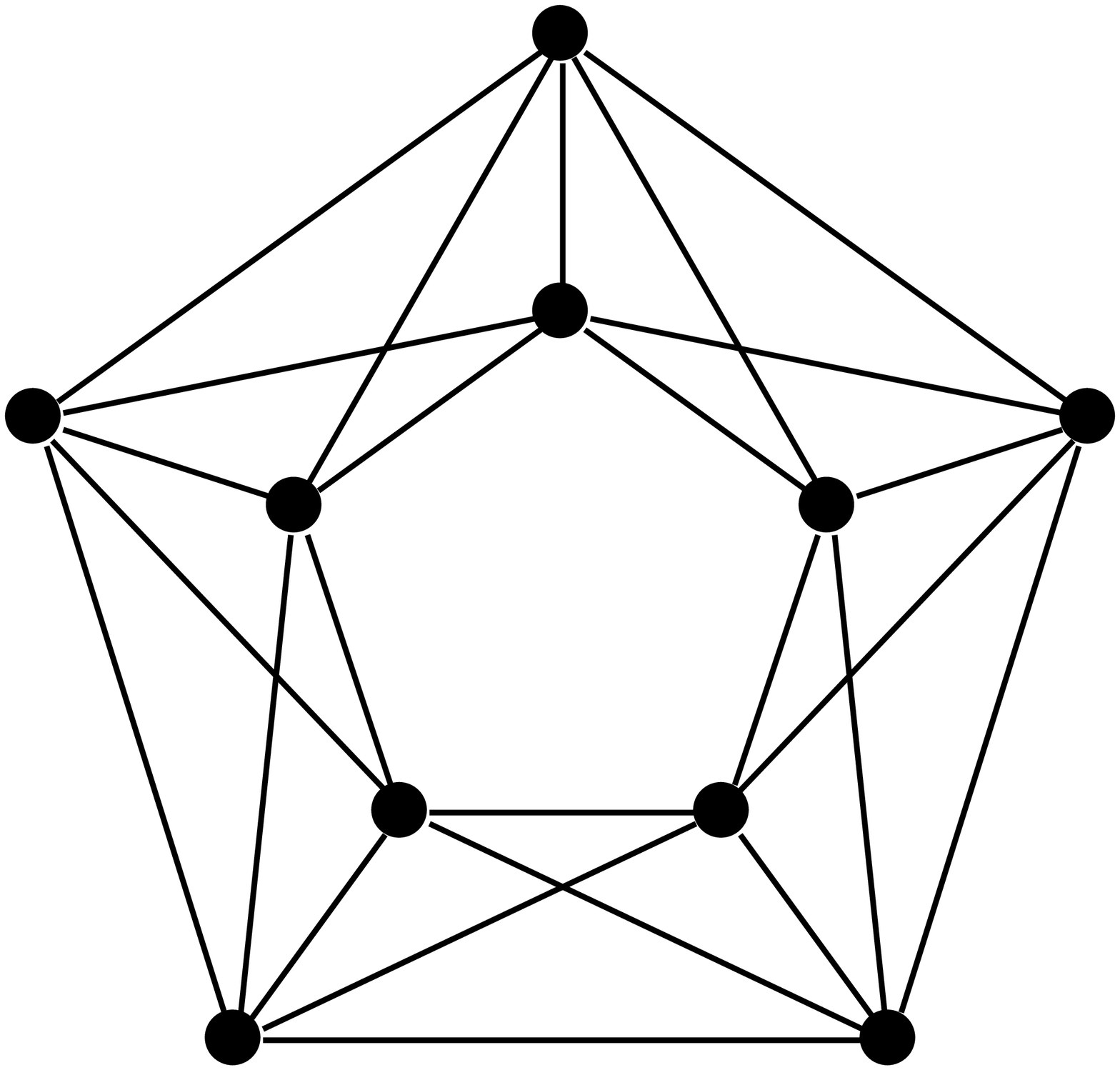,
 width=0.24\textwidth} }
 \centerline{$C_8^2$  \hspace{6cm} $C_5 \boxtimes K_2$}
\caption{The graph $C_8^2$ and $C_5 \boxtimes K_2$ .}
\end{figure}

Actually, Theorem \ref{main} is a corollary of  the following stronger
result.


\begin{theorem}\label{submain}
Assume that a connected graph $G$ is neither $C^2_8$ nor
$C_5\boxtimes K_2$. If $\Delta(G)\geq 4$ and $\omega(G)\leq
\Delta(G)-1$, then we have
$$\chi_f(G)\leq \Delta(G) -\frac{2}{67}.$$
\end{theorem}

\noindent {\bf Remark:} In the case  $\Delta(G)=3$, Heckman and
Thomas \cite{ht} conjectured that $\chi_f(G) \leq 14/5$ if $G$ is
triangle-free. Hatami and Zhu \cite{hz} proved $\chi_f(G) \leq 3-
\frac{3}{64}$ for any triangle-free graph $G$ with $\Delta(G) \leq
3$. The second and third authors showed an improved
result $\chi_f(G) \leq 3 - \frac{3}{43}$ in the previous paper \cite{lup}. Thus we need only consider
the cases $\Delta(G) \geq 4$.
 For any
connected graph $G$ with sufficiently large $\Delta(G)$ and
$\omega(G) \leq \Delta(G)-1$, Reed's result \cite{reed1}
$\chi(G)\leq \Delta(G)-1$ implies $\chi_f(G) \leq \Delta(G)-1$.
The method introduced in \cite{hz} and strengthened in \cite{lup},
has a strong influence on this paper.  The readers are encouraged to
read these two papers \cite{hz, lup}.

Let $f(k)=\inf_G\{\Delta(G)-\chi_f(G)\}$, where the infimum is taken
over all connected graphs $G$ with $\Delta(G)=k$ and not one of the
graphs listed in Theorem \ref{main}.  Since  $\chi_f(G)\geq \omega(G)$,
by taking a graph with $\omega(G)=\Delta(G)-1$, we have $f(k)\leq 1$.
Theorem \ref{submain} says
$f(k)\geq \frac{2}{67}$ for any $k\geq 4$.
 Reed's result \cite{reed1}
implies $f(k)= 1$ for sufficiently large $k$. Heckman and Thomas
\cite{ht} conjectured $f(3)=1/5$. It is an interesting problem to
determine the value of $f(k)$ for small $k$. Here we conjecture
$f(4)=f(5)=\frac{1}{3}$.  If Borodin  and Kostochka's conjecture is
true, then $f(k)= 1$ for $k\geq 9$.

 Theorem 2 is proved by induction on $k$.
Because the proof is quite long, we split the  proof into
the following two lemmas.
\begin{lemma} \label{f4}
We have $f(4) \geq \frac{2}{67}$.
\end{lemma}
\begin{lemma}\label{increase}
For each $k \geq 6$, we have $f(k)\geq \min \left\{f(k-1),
\frac{1}{2}\right\}.$ We also have  $f(5)\geq \min\left\{f(4),
\frac{1}{3}\right\}$.
\end{lemma}
It is easy to see  the combination of Lemma \ref{f4} and Lemma
\ref{increase} implies Theorem \ref{submain}. The idea of reduction
comes from the first author, who pointed out  $f(k)\geq \min \left\{f(k-1),
\frac{1}{2}\right\}$ for $k\geq 7$ based on  his  recent results
 \cite{king}. The second and third authors orginally proved 
$f(k)\geq \frac{C}{k^5}$ (for some $C>0$) using different method in the first version;
they also prove the reductions 
at $k=5,6$, which are much harder than
the case $k\geq 7$.  We do not know
whether a similar reduction exists for $k=4$.

The rest of this paper is organized as follows. In section 2, we
will introduce some notation and  prove Lemma \ref{increase}. In
section 3 and section 4, we will prove $f(4) \geq \frac{2}{67}$.

\section{Proof of Lemma \ref{increase}}

In this paper, we use the following notation. Let $G$ be a simple
graph with vertex set $V(G)$ and edge set $E(G)$. The {\it
neighborhood} of a vertex $v$ in $G$, denoted by $\Gamma_G(v)$, is
the set $\{u \colon uv \in E(G)\}$. The {\it degree} $d_G(v)$ of $v$
is the value of $|\Gamma_G(v)|$. The {\it independent
  set} (or {\it stable set}) is a set $S$ such that no edge with both
ends in $S$.  The {\it independence number} $\alpha(G)$ is the largest
size of $S$ among all the independent sets $S$ in $G$.  When $T
\subset V(G)$, we use $\alpha_G(T)$ to denote the independence number
of the induced subgraph of $G$ on $T$. Let $\Delta(G)$ be the maximum
degree of $G$.  For any two vertex-sets $S$ and $T$, we define
$E_G(S,T)$ as $\{uv \in E(G): u \in S \ \textrm{and} \ v \in T\}$.
Whenever $G$ is clear under context, we will drop the subscript $G$
for simplicity.

If $S$ is a subset of vertices in $G$, then {\it contracting} $S$
means replacing vertices in $S$ by a single fat vertex, denoted by $\underline{S}$,
 whose incident
edges are all edges that were incident to at least one vertex in
$S$, except edges with both ends in $S$. The new graph obtained by
contracting $S$ is denoted by $G/S$.  This operation is also known
as {\it identifying vertices of  $S$} in the literature. For
completeness, we allow $S$ to be a single vertex or even the empty
set. If $S$ only consists of a single vertex, then $G/S=G$; if
$S=\emptyset$, then $G/S$ is the union of $G$ and an isolated
vertex.  When $S$ consists of $2$ or $3$ vertices,  for convenience,
we write $G/uv$ for $G/\{u,v\}$ and $G/uvw$ for $G/\{u,v,w\}$; the
fat vertex will be denoted by $\underline{uv}$ and
$\underline{uvw}$, respectively. Given two disjoint subsets $S_1$
and $S_2$, we can contract $S_1$ and $S_2$ sequentially. The order
of contractions does not matter; let $G/S_1/S_2$ be the resulted
graph. We use $G-S$ to denote the subgraph of $G$ induced by
$V(G)-S$.

In order to prove Lemma \ref{increase}, we need use the following theorems
due to King \cite{king}.

\begin{theorem}[King \cite{king}]
\label{king0}
If a graph $G$ satisfies $\omega(G) > \frac{2} {3} (\Delta(G) + 1)$,
then $G$ contains a stable set $S$ meeting every maximum clique.
\end{theorem}
\begin{theorem}[King \cite{king}]
\label{king}
For a positive integer $k$, let $G$ be a graph with vertices
partitioned into cliques $V_1,\ldots,V_r$. If for every $i$ and
every $v \in  V_i$, $v$ has at most $\min\{k, |V_i|-k\}$ neighbors
outside $V_i$, then $G$ contains a stable set of size $r$.
\end{theorem}

\begin{lemma} \label{6to5}
Suppose that  $G$ is a connected graph with $\Delta(G)\leq 6$ and $\omega(G)\leq 5$.
Then there exists an independent set  meeting  all induced copies
of  $K_5$ and $C_5\boxtimes K_2$.
\end{lemma}

\noindent
{\bf Proof:} We first show that there exists an independent set
meeting all copies of $K_5$. If $G$ contains no $K_5$, then this is
trivial. Otherwise, we can apply Theorem \ref{king0} to get the
desired independent set since $\omega(G)>\frac{2}{3}(\Delta(G)+1)$
is satisfied.

Now we prove the Lemma by contradiction. Suppose the Lemma is false.
Let $G$ be a  minimum counterexample (with the smallest number of
vertices).  For any independent set $I$, let $C(I)$ be the number of
induced copies of $C_5\boxtimes K_2$ in $G-I$. Among all independent
sets which  meet all copies of $K_5$, there exists one such
independent set $I$ such that $C(I)$ is minimized.

Since $C(I)>0$, there is  an induced copy of $C_5\boxtimes K_2$ in
$G-I$; we use $H$ to denote it.  In $C_5\boxtimes K_2$, there is a
unique perfect matching such that identifying
 the two ends of each edge in this matching results  a $C_5$.
An edge in this unique matching is called  a {\em canonical } edge.
We define a new graph $G'$ as follows: First we  contract all
canonical edges in $H$ to get a $C_5$, where its vertices are called {\em fat}
vertices.
Second we add five edges
turning the $C_5$ into a $K_5$. Observe that  each vertex in this
$C_5$ can have at most two neighbors in $G-H$ and  $\Delta(G') \leq
6$. We will consider the following four cases.

\noindent{\bf Case 1:} There is a $K_6$ in the new graph $G'$. Since
the original graph $G$ is $K_6$-free, the $K_6$ is formed by the
following two possible ways.

\begin{figure}[htbp]
\centerline{ {\psfig{figure=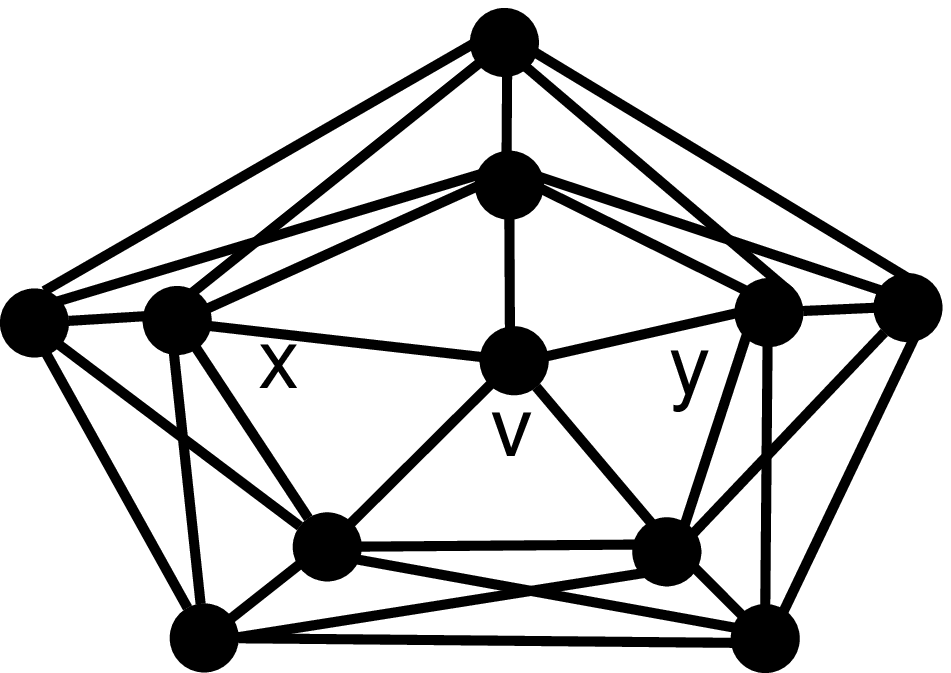, width=0.32\textwidth}}
\hfil \psfig{figure=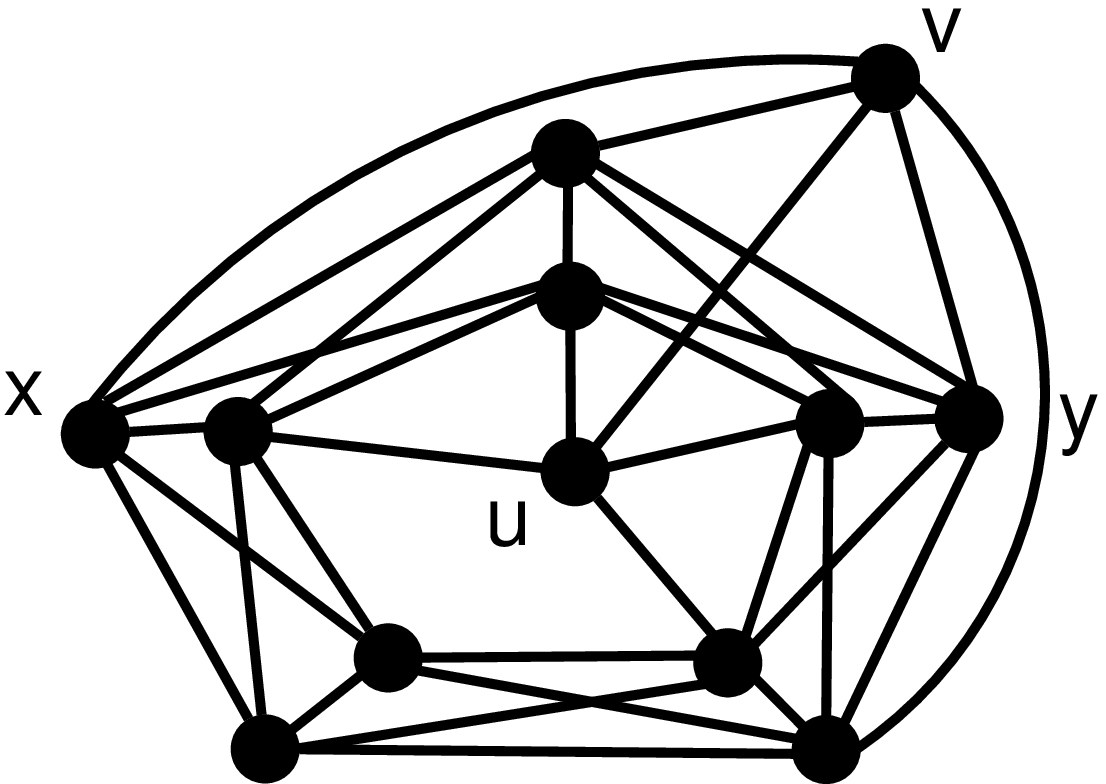, width=0.32\textwidth}}
  \begin{multicols}{2}
\caption{Subcase 1a.} 
\label{fig:1a} 
\newpage
\caption{Subcase 1b.} 
\label{fig:1b} 
\end{multicols}
\end{figure}

{\bf Subcase 1a:} This $K_6$ contains $5$ fat vertices. By the symmetry of $H$,
there is an induced  $C_5$ in $H$
 such that the vertices in $C_5$ contain a common neighbor vertex $v$ in $G \setminus V(H)$,  see Figure \ref{fig:1a}.
 Since $H$ is $K_5$-free, we can find $x, y$  in this $C_5$ such that $x,y$ is a non-edge.
Let $I':=(I\setminus\{v\})\cup\{x, y\}$;  $I'$ is also an
independent set.  Observe that  $v$  is not  in any $K_5$ in $G-I'$. Thus the set $I'$ is also an independent set and meets every
$K_5$ in $G$. Since $C_5\boxtimes K_2$ is a $5$-regular graph, any
copy of $C_5\boxtimes K_2$ containing $v$ must contain at least one
of $x$ and  $y$. Thus, $C(I')<C(I)$. Contradiction!

{\bf Subcase 1b:}  This $K_6$ contains $4$ fat vertices. 
Let $u,v$ be the other two vertices.
By the  symmetry of $H$, there is a unique way to connect $u$ and $v$ to $H$
as shown by Figure \ref{fig:1b}.
Since $uv$ is an edge, one of $u$ and $v$ is not in $I$. We assume
$u \not \in I$. Let $\{x,y\} \subset \Gamma_G(v) \cap V(H)$ as shown
in Figure \ref{fig:1b} and $I'=I \setminus \{v\} \cup \{x,y\} $.
Observe that $I'$ is an independent set and $v$ is not in a $K_5$ in $G-I'$.  Thus $I'$ is an independent set meeting each $K_5$ in $G$.
Since each $C_5 \boxtimes K_2$ containing $v$ must contain one of
$x$ and $y$. Thus $C(I')<C(I)$. Contradiction!

\begin{figure}[htbp]
\centerline{ {\psfig{figure=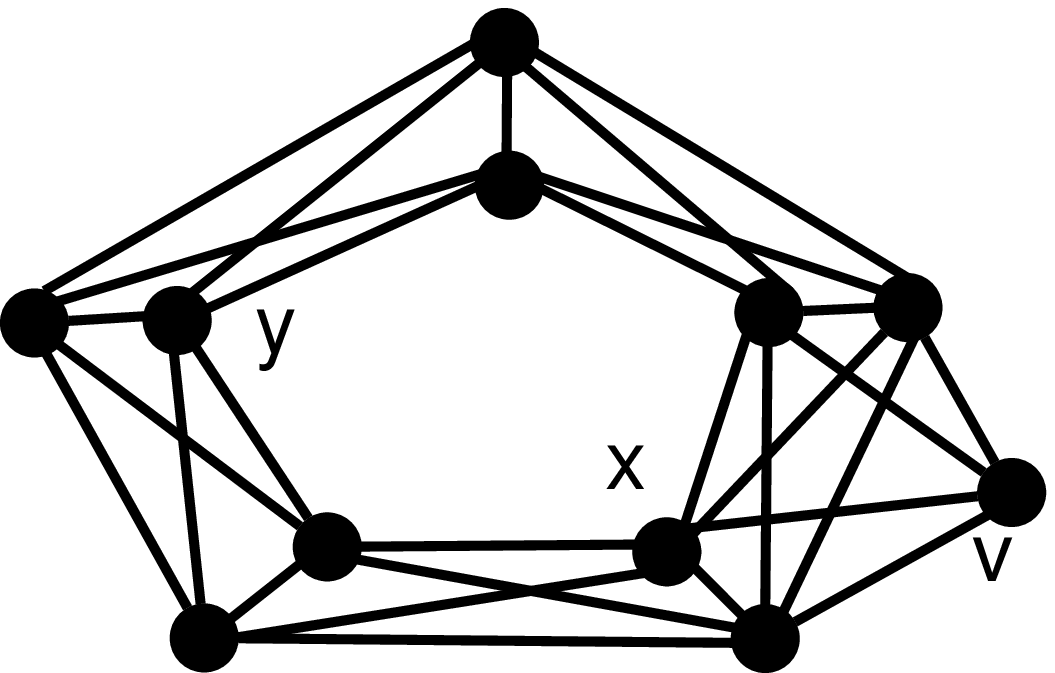, width=0.32\textwidth}}
\hfil \psfig{figure=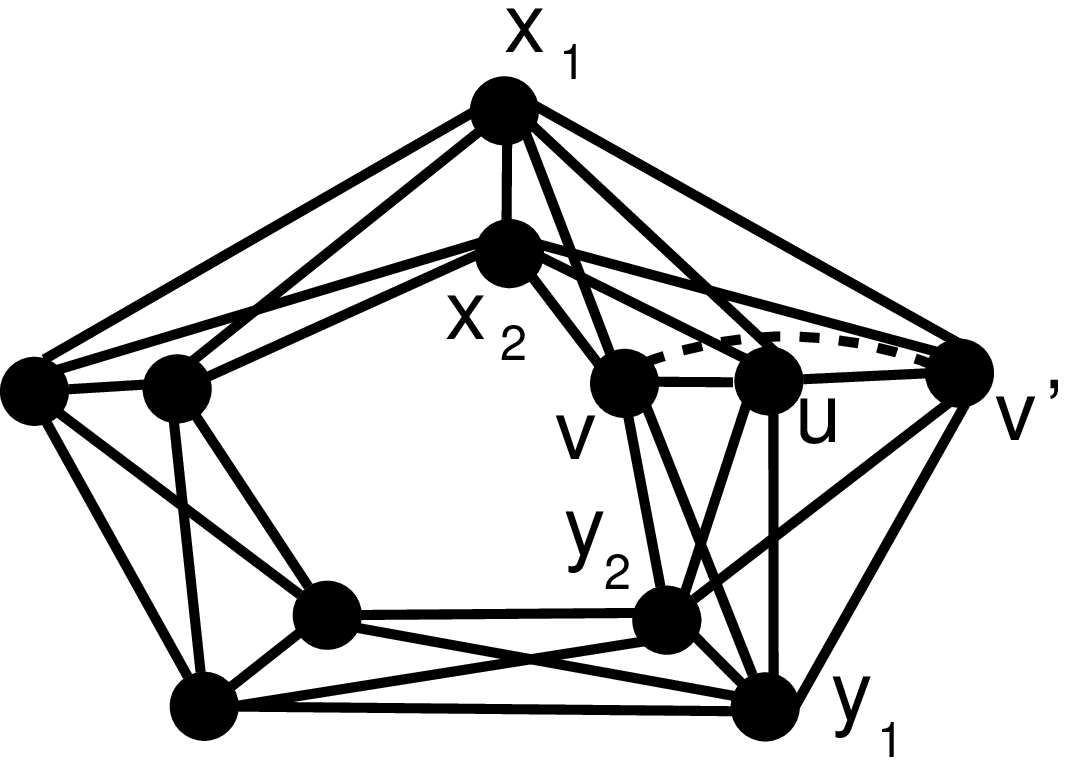, width=0.32\textwidth}}
  \begin{multicols}{2}
\caption{Case 2.} 
\label{fig:2} 
\newpage
\caption{Case 3.} 
\label{fig:3} 
\end{multicols}
\end{figure}

\noindent {\bf Case 2:}
 There is a $K_5$   intersecting $H$ with $4$
vertices. Let $v$ be the vertex of this $K_5$ but not in $H$, see Figure \ref{fig:2}.
We have two subcases.

{\bf Subcase 2a:} The vertex $v$ has another neighbor $y$ in $H$ but
not in this $K_5$. Since $H$ is $K_5$-free, we can select a vertex
$x$ in this $K_5$ such that $xy$ is not an edge of $G$. Let
$I':=I\setminus \{v\}\cup \{x,y\}$. Note that $v$ is not in a $K_5$ in $G-I'$, and $I'$ is an independent set. Thus $I'$ is an
independent set meeting each $K_5$ in $ G$.  Since  any  $C_5
\boxtimes K_2$ containing $v$ must contain one of $x$ and $y$, we
have  $C(I')<C(I)$. Contradiction!

{\bf Subcase 2b:} All neighbors of  $v$ in $H$ are in this $K_5$.
Let $x$ be any vertex in this $K_5$ other than $v$, and
$I':=I\setminus \{v\}\cup \{x\}$. In this case, there is only one
$K_5$ containing $v$. Thus, $I'$ is also an independent set meeting
every copy of $K_5$ in $G$. Observe that $\Gamma_G(v) \setminus
\{x\}$ is disconnected. If $v \in H'=C_5 \boxtimes K_2$, then
$\Gamma_G(v) \cap H'$ is connected.  Thus $v$ is not in a $C_5\boxtimes
K_2$ in $G-I'$ and $C(I')< C(I)$. Contradiction!

\noindent {\bf Case 3:}  There is an induced subgraph $H'$
isomorphic to $C_5\boxtimes K_2$ such that  $H'$ and $H$ are
intersecting, see Figure \ref{fig:3}.
Since $V(H)\cap V(H')\not =\emptyset$ and $H\not=H'$, we can find a
canonical edge $uv$ of $H$ and a canonical edge $uv'$ of $H'$ such
that $v\not\in V(H')$ and $v'\not\in V(H)$.
 If $vv'$ is a non-edge, then let $I':=I\setminus \{v'\} \cup \{u\}$.
It is easy to check $I'$ is still an independent set. We also
observe that any possible $K_5$ containing $v'$ must also contain
$u$. Thus, $I'$ meets every copy of $K_5$ in $G$.  We have $v'$ in no $C_5 \boxtimes K_2$ in $G-I'$ since $vv'$ is not an edge.   We therefore
get $C(I')<C(I)$. Contradiction! If $vv'$ is an edge, then locally
there are two $K_5$ intersecting at $u$, $v$, and $v'$;  say the
other four vertices are $x_1, x_2, y_1, y_2$, where two cliques are
$\{x_1,x_2,u,v,v'\}$ and $\{y_1,y_2,u,v,v'\}$, see Figure
\ref{fig:3}. Let $I'=I\cup\{x_1,y_1\}\setminus \{v'\}$. Note that
$I'$ is an independent set and $v'$ is not in a $K_5$ in $G-I'$. Thus
$I'$ is an independent set meeting each $K_5$ in $G$. Observe that
 any copies of $C_5\boxtimes K_2$  containing $v'$ must  contain one
 of $x_1$ and $y_1$; we have  $C(I')<C(I)$. Contradiction!.

\noindent
{\bf Case 4:}  This is the remaining case, $G'$ is $K_6$-free.  We
have $\omega(G')\leq 5$ and $|V(G')| < |V(G)|$. By the minimality of
$G$, there is an independent set $I'$ of $G'$ meeting every copy of
$K_5$ and $C_5\boxtimes K_2$. In $I'$, there is a unique vertex $x$
of the $K_5$ obtained from contracting  canonical edges of $H$. Let
$uv$ be the canonical edge corresponding to $x$. Let $I''=I'
\setminus \{x\} \cup \{u\}$, we get an independent set $I''$ of $G$.
Note that any $v \in H \setminus \{u\}$ is not in any $K_5$ of
$G-I''$ by Case 2 as well as not in any $C_5 \boxtimes K_2$ of
$G-I''$ by Case 3. Thus $I''$ hits each $K_5$ in $G$ and $C(I'')=0$.
Contradiction! \hfill $\square$

The following lemma extends Theorem \ref{king0} when $\omega(G)=4$; a similar result was proved independently in \cite{cek}.

\begin{lemma} \label{5to41}
Let $G$ be a connected graph with $\Delta(G)\leq 5$ and $\omega(G)\leq 4$.
If $G \not = C_{2l+1} \boxtimes K_2$ for some $l \geq 2$,
then there is an independent set $I$ hitting all copies of $K_4$ in $G$.
\end{lemma}
{\bf Proof:} We will prove it by contradiction. If the lemma is
false, then let $G$ be a minimum counterexample. If $G$ is
$K_4$-free, then there is nothing to prove. Otherwise, we consider
the clique graph $\C(G)$, whose edge set is the set of  all edges
appearing in some copy of $K_4$. Because of $\Delta(G)=5$, here are
all possible connected component of $\C(G)$.

\begin{enumerate}
\item $C_{t}\boxtimes K_2$ for $t\geq 4$.  If this type occurs, then every vertex
in $C_{t}\boxtimes K_2$ has degree $5$; thus, this is the entire
graph $G$. If $t$ is even, then we can find an independent set $I$
meeting every $K_4$. If $t$ is odd, then it is impossible to find
such an independent set. However, this graph is excluded from the
assumption of the Lemma.

\item $P_t\boxtimes K_2$ for $t\geq 3$. In this case, all internal vertices
have degree $5$ while the four end vertices  have degree $4$.
Consider a new graph $G'$ which is obtained by deleting all internal
vertices and adding four edges to make the four end vertices as a
$K_4$. It is easy to check $\Delta(G')\leq 5$ and $\omega(G')\leq
4$. Since $|G'|<|G|$, there is an independent set $I$ of $G'$
meeting every copy of $K_4$ in $G'$. Note  that there is exactly one
end vertex in $I$. Observe that any one end vertex  can be extended
into a maximal independent set meeting every copy of $K_4$ in
$P_t\boxtimes K_2$. Thus, we can extend $I$ to an independent set
$I'$ of $G$ such that $I'$ meets every copy of $K_4$ in $G$. Hence,
this type of component does not occur in $\C(G)$.

\item There are four other types listed in Figure \ref{k4}.
\begin{figure}[htbp]
\centerline{ {\psfig{figure=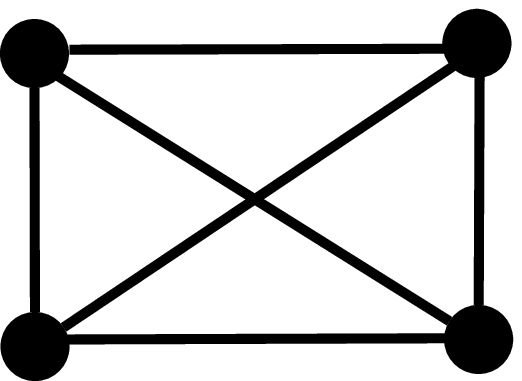, width=0.2\textwidth}}
\hfil
\psfig{figure=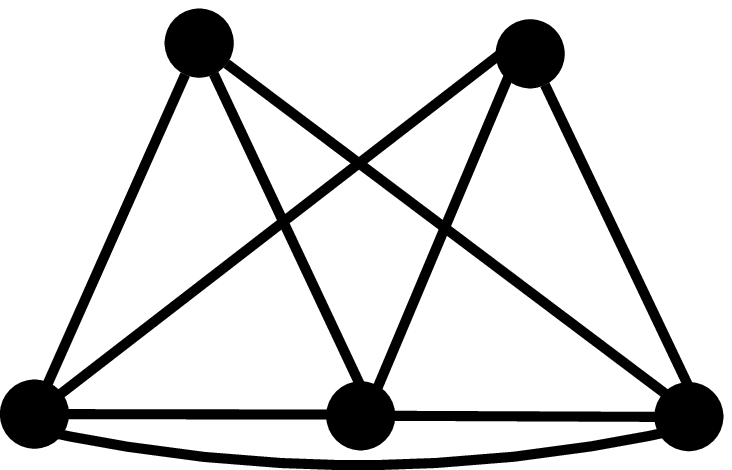, width=0.2\textwidth}\hfil
   {\psfig{figure=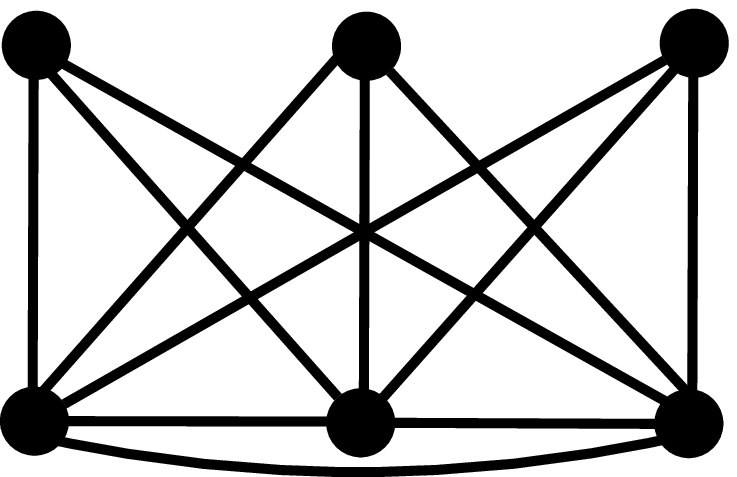, width=0.2\textwidth}}
    \hfil \psfig{figure=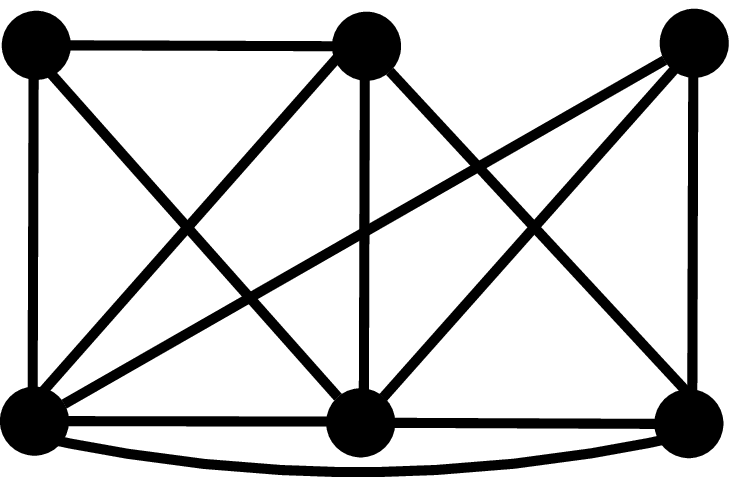, width=0.2\textwidth}}
\caption{All types of components in the clique graph $C(G)$.}
 \label{k4}
\end{figure}

For each  component $C_i$ in $\C(G)$, let $V_i$ be the set of common
vertices in all $K_4$'s of $C_i$; for the leftmost figure in Figure
\ref{k4}, $V_i$ is the set of all 4 vertices; for the  middle two
figures, $V_i$ is the set of bottom three vertices; for the
rightmost figure, $V_i$ consists of the left-bottom vertex and the
middle-bottom vertex. Note that all $V_i$'s are pairwise disjoint.
Let $G'$ be the induce subgraph of $G$ on $\cup_i V_i$. Note that
$G'$ does not contains any vertex in $C_i\setminus V_i$. By checking
each type, we find out that
 for each $i$ and each $v\in V_i$, $v$ has at most $\min\{2,|V_i|-2\}$
neighbors outside $V_i$ in $G'$  (not in $G$!). Applying Theorem
\ref{king} to $G'$, we conclude that there exists an independent set
$I$ of $G'$ meeting every $V_i$; thus $I$ meets every $K_4$ in $G$.
Contradiction!
\end{enumerate}
\hfill $\square$

\begin{lemma} \label{5to42}
Let $G$ be a connected graph with $\Delta(G)\leq 5$ and
$\omega(G)\leq 4$. If $G \not = C_{2l+1} \boxtimes K_2$ for some $l
\geq 2$, then there exists an independent set  meeting  all induced
copies of  $K_4$ and $C_8^2$.
\end{lemma}

\noindent
{\bf Proof:} We will use proof by contradiction.  Suppose the Lemma
is false. Let $G$ be a minimum counterexample (with the smallest
number of vertices).  For any independent set $I$, let $C(I)$ be the
number of induced copies of $C_8^2$ in $G-I$.  Among all independent
sets which meet all copies of $K_4$, there exists an independent set
$I$ such that $C(I)$ is minimized. Since $C(I)>0$, let $H$ be a copy
of $C_8^2$ in $G-I$. The vertices of $H$ are listed by $u_i$ for
$i\in {\mathbb Z}_8$  anticlockwise such that $u_iu_j$ is an edge of
$H$ if and only if $|i-j|  \leq 2$. The vertex $v_{i+4}$ is the {\em
antipode} of $v_i$ for any $i\in {\mathbb Z}_8$.

\noindent
{\bf Case 1:} There exists a vertex $v\not\in V(H)$ such that $v$ has
 five neighbors in $H$. By the Pigeonhole Principle, $\Gamma(v)$ contains
a pair of antipodes. Without loss of generality, say $u_0, u_4\in
\Gamma(v)$. If the other three neighbors of $v$ do not form a
triangle, then we let $I':=I\setminus \{v\} \cup \{u_0,u_4\}$; note
that $v$ is not in any $K_4$ of $G-I'$. Thus $I'$ is an independent
set meeting every copy of $K_4$. Since every copy of $C_8^2$
containing $v$ must contain one of $u_0$ and $u_4$, we have
$C(I')<C(I)$. Contradiction! Hence, the other three neighbors of $v$
must  form a triangle. Without loss of generality, we can assume
that the  three neighbors are $u_1$, $u_2$, and $u_3$. Now we let
$I'':=I\setminus \{v\} \cup \{u_0,u_3\}$; note that $v \not \in K_4
\subset G-I'$. Thus  $I''$ is also an independent set meeting every
copy of $K_4$ of $G$. Since every copy of $C_8^2$ containing $v$
must contain one of $u_0$ and $u_3$, we have $C(I'')<C(I)$.
Contradiction!

\noindent
{\bf Case 2:} There exists a vertex $v\not\in V(H)$ such that $v$
has exactly four neighbors in $H$. Since $H$ is $K_4$-free, we can
find $u_i, u_j\in \Gamma(v)\cap V(H)$ such that $u_iu_j$ is a
non-edge. Let  $I':=I\setminus \{v\} \cup \{u_i,u_j\}$; $I'$ is also
an independent set.   Note that $\Gamma(v)\setminus \{u_i, u_j\}$
can not be a triangle, $v$ is not in any $K_4 \subset G-I'$. Thus
$I'$ meets every copy of $K_4$. Since every copy of $C_8^2$
containing $v$ must contain one of $u_i$ and $u_j$, we have
$C(I')<C(I)$. Contradiction!

\noindent
{\bf Case 3:} There exists a vertex $v\not\in V(H)$ such that $v$
has exactly three neighbors in $H$. If the $3$ neighbors do not form
a triangle, then choose $u_i, u_j\in \Gamma(v)\cap V(H)$ such that
$u_iu_j$ is a non-edge. Note that $\Gamma(v)\setminus \{u_i, u_j\}$
can not be a triangle; $v$ is not in any $K_4 \subset G-I'$.  Let
$I':=I\setminus \{v\} \cup \{u_i,u_j\}$; $I'$ is also an independent
set meeting every copy of $K_4$. Since every copy of $C_8^2$
containing $v$ must contain one of $u_i$ and $u_j$, we have
$C(I')<C(I)$. Contradiction! Else, the three neighbors form a
triangle; let $u_i$ be one of them and $I':=I\setminus \{v\} \cup
\{u_i\}$;  $v$ is not in any $K_4 \subset G-I'$. Thus $I'$ is  an
independent set meeting every copy of $K_4$.   Note that
$\Gamma(v)\setminus \{u_i\}$ has only two vertices in $H$. The
induced graph on $\Gamma(v)\setminus \{u_i\}$ is disconnected.
However, for any vertex $v$ in $H'=C_8^2$,  the subgraph induced by
$\Gamma_G(v) \cap H'$ is a $P_4$. There is no $C_8^2$ in $G-I'$
containing $v$. Thus, $C(I')<C(I)$. Contradiction!

\noindent
{\bf Case 4:} Every vertex outside $H$ can have at most $2$
neighbors in $H$. We identify each pair of antipodes of $H$ to get a
new graph $G'$ from $G$. After identifying, $H$ is turned into a
$K_4$; where the vertices of this $K_4$ are referred as fat
vertices.

\begin{figure}[htbp]
\centerline{  \psfig{figure=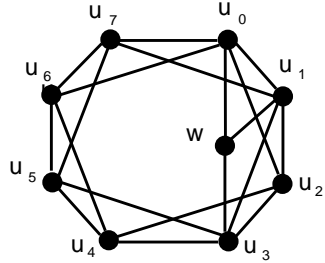, width=0.3\textwidth}}
 \caption{Subcase 4a.}
 \label{subcase4a}
\end{figure}

 {\bf Subcase 4a:} $G'\not=C_{2l+1}\boxtimes K_2$.
 Observe  $\Delta(G')\leq 5$. We claim $G'$ is $K_5$-free.
Suppose not. Since every vertex in $H$ has at most one neighbor
outside $H$, then each fat vertex can have at most two neighbors
outside $H$.  Recall that the original graph $G$ is $K_5$-free. If
$G'$ has some $K_5$, then this $K_5$ contains either $3$ or $4$ fat 
vertices. Let $w$ be one of the other vertices in this $K_5$. 
We get $w$ has at least three neighbors in
$H$. However, this is covered by Case 1, Case 2, or Case 3. Thus,
$G'$ is $K_5$-free.
 Since $|G'|<|G|$, by the minimality of $G$, $G'$ has an independent set $I'$
meeting every copy of $K_4$ and $C_8^2$ in $G'$.  There is exactly
one fat vertex in $I'$. Now replacing this fat vertex by its
corresponding pair of antipodal vertices, we get an independent set
$I''$; we assume the pair of antipodal vertices are $u_2$ and $u_6$.
It is easy to check that $I''$ is an independent set of $G$. Next we
claim any $v \in V(H) \setminus \{u_2,u_6\}$ is neither in a $K_4
\subset V(G)-I''$ nor in a $C_8^2 \subset V(G)-I''$. Suppose there is
some $v$ such that $v \in K_4 \subset G-I''$. Recall each $v \in
V(H)$ has at most one neighbor outside $H$ and $H$ is $K_4$-free;
there is some $w \not \in V(H)$ such that $w$ has at least three
neighbors in $H$. This is already considered by Case 1, Case 2, or
Case 3.  We are left to show that $v \not \in C_8^2 \subset G-I''$
for each $v \in V(H) \setminus \{u_2,u_6\}$.  If not, there exists 
a copy $H'$ of  $C_8^2$ in $G-I''$ containing $v$. Note $H'$ is $4$-regular,
any vertex in $H'$ can have at most one neighbor in $I''$; in particular,
$v\not=u_0, u_4$. Without loss of generality, we assume $v=u_3$.
 Then there is a vertex $w \not \in V(H)$ such
that $u_3w$ is an edge, see Figure \ref{subcase4a}. Observe that the
neighborhood of each vertex of an induced $C_8^2$ is is a $P_4$.
Since $u_1u_4$ and $u_1u_5$ are two non-edges, we have $wu_1$ being
an edge.  Observe $\Gamma_G(u_1)=\{u_7,u_0,u_2,u_3,w\}$. Since $u_2
\not \in H'$, we have $u_0 \in H'$; $u_0$ has two neighbors ($u_2$ and $u_6$)
outside $H'$, contradiction!
  Therefore, $I''$ meets every copy of $K_4$ and $C_8^2$ in $G$.
 Contradiction!

{\bf Subcase 4b:} $G'=C_{2l+1}\boxtimes K_2$. The graph $G$ can be
recovered from $G'$. It consists of an induced subgraph $H=C_8^2$
and an induced subgraph $P_{2l-1}\boxtimes K_2$. For each vertex $u$
in $H$, there is exactly one edge connecting it to one of the four
end vertices of $P_{2l-1}\boxtimes K_2$; for each end vertex $v$ of
$P_{2l-1}\boxtimes K_2$, there are exactly two edges connecting $v$
to the vertices in $H$. First, we take any maximum independent set
$I'$ of $P_{2l-1}\boxtimes K_2$. Observe that $I'$ has exactly two
end points of $P_{2l-1} \boxtimes K_2$; so $I'$  has exactly four
neighbors in $H$. In the remaining four vertices of $H$, there
exists a non-edge $u_iu_j$ since $H$ is $K_4$-free. Let
$I:=I'\cup\{u_i,u_j\}$. Clearly $I$ is an independent set of $G$
meeting every copy of $K_4$ and $C_8^2$. Contradiction! \hfill
$\square$

We are ready to prove  Lemma \ref{increase}.

\noindent
{\bf Proof of Lemma \ref{increase}:} We need  prove
for $k\geq 5$ and any connected graph $G$ with $\Delta(G)=k$
and $\omega(G)\leq k-1$ satisfies
\begin{equation}
  \label{eq:rec}
\chi_f(G)\leq k- \min\left\{f(k-1),\frac{1}{2}\right\}.
\end{equation}
 If $\omega(G)\leq k-2$,
then by inequality (\ref{eq:1}), we have
$$\chi_f(G)\leq \frac{\Delta(G)+\omega(G)+1}{2}\leq k-\frac{1}{2}.$$
Thus, inequality (\ref{eq:rec}) is satisfied. From now on, we assume
$\omega(G)=\Delta(G)-1$.

For $\Delta(G)=k\geq 6$ and $\omega(G)=k-1$, the condition
$\omega(G)>\frac{2}{3}(\Delta(G)+1)$ is satisfied. By Theorem
\ref{king0}, $G$ contains an independent set meeting every maximum
clique. Extend this independent set to a maximal independent set and
denote it by $I$.
 Note that $\Delta(G-I) \leq k-1$ and $\omega(G-I)\leq k-2$.

\noindent {\bf Case 1:} $k\geq 7$. From the definition of $f(k-1)$,
we have $\chi_f(G-I) \leq \Delta(G-I)-f(k-1)$. Thus,
$$
\chi_f(G) \leq \chi_f(G-I)+1 \leq k-1-f(k-1)+1=k-f(k-1).
$$
Thus, we have $f(k) \geq \min\{f(k-1),1/2\}$.

\noindent {\bf Case 2:} $k=6$. By Lemma \ref{6to5}, we can find an
independent set   meets every copy of $K_5$ and $C_5\boxtimes K_2$;
we extend this independent set  as a maximal independent set $I$.
Note that
 $G-I$ contains no induced subgraph isomorphic $C_5\boxtimes K_2$.
We have $\chi_f(G-I) \leq 5-f(5)$;  it implies $\chi_f(G)\leq
6-f(5)$.
 Thus, $f(6) \geq \min\{f(5),1/2\}$ and we are done.

\noindent
{\bf Case 3:} $k=5$. If
$G=C_{2l+1}\boxtimes K_2$ for some $l\geq 3$; then $G$ is
vertex-transitive and $\alpha(G)=l$. It implies that
$$\chi_f(G)=\frac{|V(G)|}{\alpha(G)}=4+\frac{2}{l}\leq 5-\frac{1}{3}.$$

If $G\not =C_{2l+1}\boxtimes K_2$, then by Lemma \ref{5to42}, we can
find an independent set   meeting every copy of $K_4$ and $C_8^2$; we
extend it as a maximal independent set  $I$. Note that
 $G-I$ contains no induced subgraph isomorphic $C_8^2$.
We have $\chi_f(G-I) \leq 4-f(4)$; it implies $\chi_f(G)\leq 5-f(4)$.
 Thus, $f(3) \geq \min\{f(4),1/3\}$ and we are finished.
\hfill $\square$

\section{The case $\Delta(G)=4$}

To prove $f(4)\geq \frac{2}{67}$, we will use an approach which is
similar to those in \cite{hz, lup}.  We will construct 133 4-colorable auxiliary graphs, and from these colorings we will construct a 134-fold coloring of $G$ using 532 colors.

It suffices to prove that the
minimum counterexample does not exist.

Let $G$ be a graph with the smallest number of vertices and
satisfying
\begin{enumerate}
\item $\Delta(G)=4$ and $\omega(G)\leq 3$;
\item $\chi_f(G)> 4-\frac{2}{67}$;
\item $G\not=C_8^2$.
\end{enumerate}
%

By the minimality of $G$, each vertex in $G$ has degree either $4$
or $3$.
 To prove Lemma \ref{f4},
 we
will show $\chi_f (G) \leq 4 - \frac{2}{67}$, which gives us the
desired contradiction.

For a given vertex $x$ in $V(G)$, it is easy to color its
neighborhood $\Gamma_G(x)$ using $2$ colors. If $d_G(x)=3$, then we
pick a non-edge $S$ from $\Gamma_G(x)$ and color the two vertices in
$S$ using  color 1. If $d_G(x)=4$ and $\alpha(\Gamma_G(x))\geq 3$,
then we pick an independent set $S$ in $\Gamma_G(x)$  of size 3 and
assign the  color 1 to each vertex in $S$. If $d_G(x)=4$ and
$\alpha(\Gamma_G(x))=2$, then we pick two
 disjoint non-edges $S_1$ and $S_2$ from $\Gamma_G(x)$ ;  we assign color 1 to each vertex in $S_1$
and  color 2 to  each vertex in $S_2$.

The following Lemma shows that $G$ has a key property, which
eventually implies that this local coloring scheme  works
simultaneously for  $x$ in a large subset of   $V(G)$.

\begin{lemma}\label{D41}
 For each $x \in V(G)$
with $d_G(x)=4$ and $\alpha(\Gamma_G(x))=2$, there exist two
vertex-disjoint non-edges $S_1(x), S_2(x) \subset \Gamma_G(x)$
satisfying the following property. If we contract $S_1(x)$ and
$S_2(x)$, then the resulting graph $G/S_1(x) /S_2(x)$ contains
neither  $K^-_{5}$ nor $G_0$. Here  $K_5^-$ is the graph obtained
from $K_5$ by removing one edge and $G_0$ is the graph shown in
Figure \ref{fig:H13}.
\begin{figure}[htbp]
  \centering
\psfig{figure=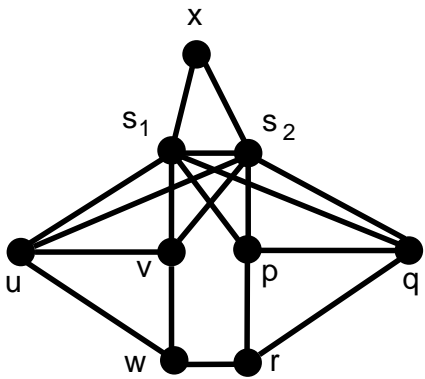,
width=0.4\textwidth}
  \caption{The graph $G_0$.}
  \label{fig:H13}
\end{figure}
\end{lemma}
The proof of this lemma is quite long and we will present its proof
in section 4.

For each vertex $x$ in $G$, we associate
 a small set of vertices  $S(x)$ selected from $\Gamma_G(x)$ as follows. If
$d_G(x)=3$,  then let $S(x)$ be the endpoints of a non-edge in
$\Gamma_G(x)$ and label the vertices in $S(x)$ as 1; if $d_G(x)=4$
and $\alpha(\Gamma_G(x))\geq 3$, then let $S(x)$ be any independent
set of size $3$ in $\Gamma_G(x)$ and label all vertices in $S(x)$ as
1; if $d_G(x)=4$ and $\alpha(\Gamma_G(x))=2$, then let
$S(x)=S_1(x)\cup S_2(x)$, where $S_1(x)$ and $S_2(x)$ are guaranteed
by Lemma \ref{D41}; we label the vertices in $S_1(x)$ as 1 and the
vertices in $S_2(x)$ as 2.   For any $x\in V(G)$, we have
$|S(x)|=2$, $3$, or $4$.

The following definitions depend on the choice of $S(\ast)$, which
is assumed to be fixed through this section.  For $v \in G$ and $j
\in \{1,2,3\}$, we define
$$
N_G^j(v) = \{u| \ \mbox{there is a path} \  vv_0 \ldots v_{j-2}u \
\mbox{in}\  G \ \mbox{of length} \   j \  \mbox{such that} \  v_0
\in S(v) \ \mbox{and} \  v_{j-2}  \in  S(u)\}.
$$

\begin{figure}[htbp]
  \centering{
\psfig{figure=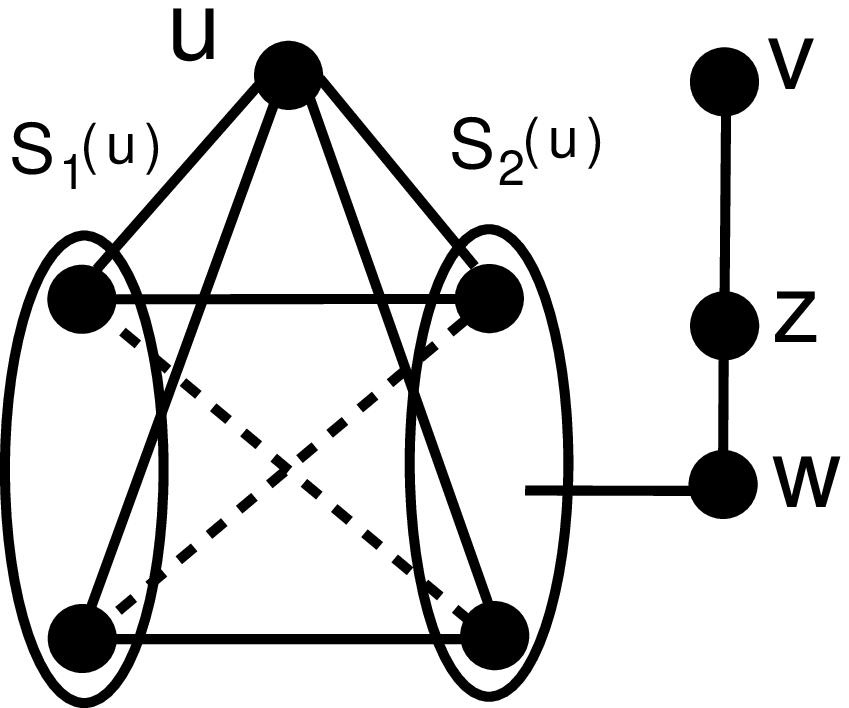, width=0.25\textwidth} \hfil
\psfig{figure=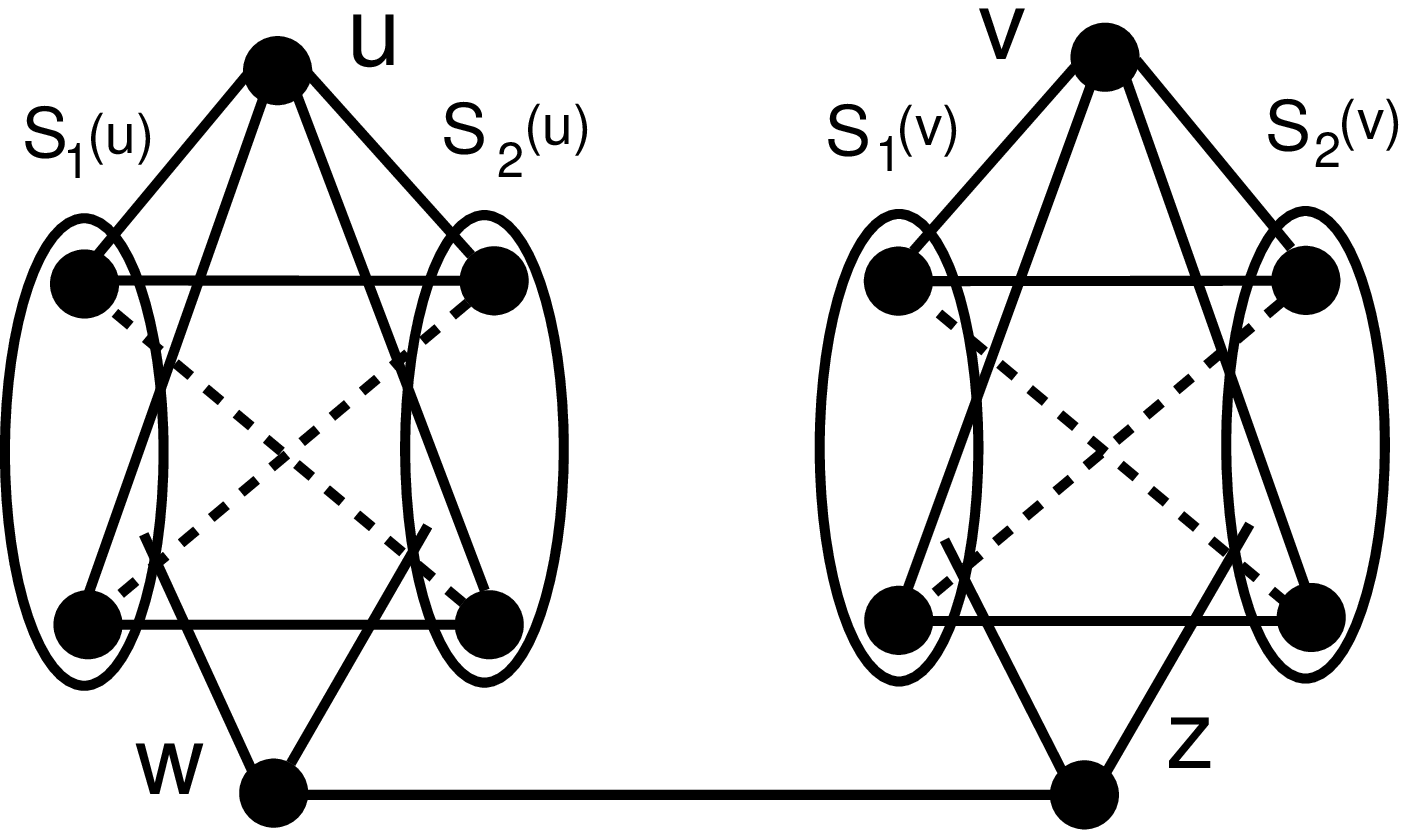, width=0.34\textwidth} \hfil
\psfig{figure=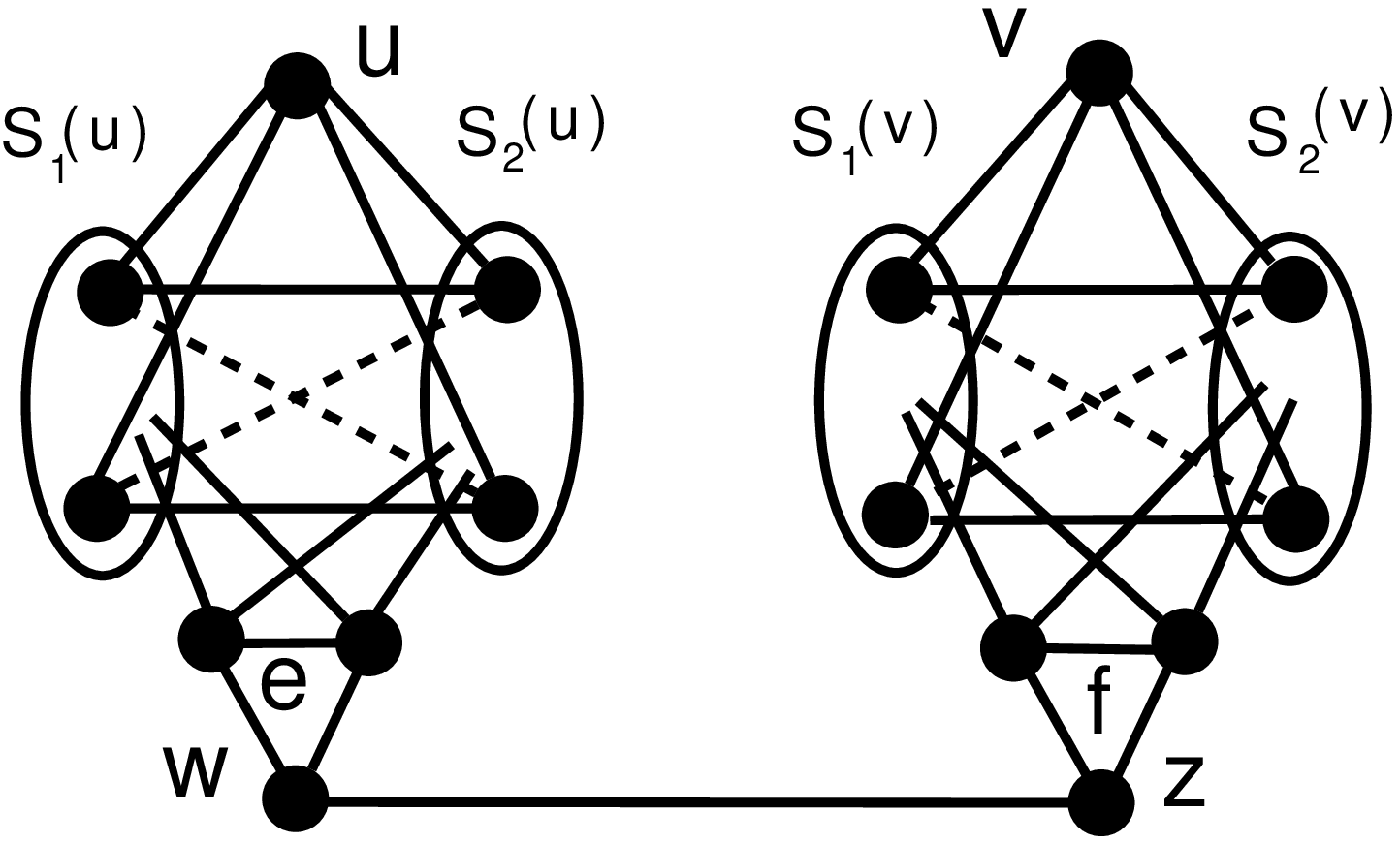, width=0.32\textwidth}}
\begin{multicols}{3}
\caption{\! 4-th neighborhood.\!\!\!} 
\label{fig:4n} 
\newpage
\caption{\! 5-th neighborhood.\!\!\!} 
\label{fig:5n} 
\newpage
\caption{\! 7-th neighborhood.\!\!\!} 
\label{fig:7n} 
\end{multicols}
\end{figure}

We now define $N_G^j(u)$ for $j\in \{4,5,7\}$; each $N_G^j(u)$ is a subset of the $j$th neighborhood of $u$.  For $j=4$,  $v \in N_G^4(u)$ if  $d_G(u)=4$,
$\alpha(\Gamma_G(u))=2$, $u$ and $v$ are connected as shown in 
Figure \ref{fig:4n}; otherwise $N_G^4(u)=\emptyset$.
In Figure \ref{fig:4n}, $w$ is connected to one of the  two vertices in $S_2(u)$.
Similarly, in  Figure \ref{fig:5n} and \ref{fig:5n}, a vertex is connected to a group
of vertices means it is connected to any vertex in this group. 

For $j=5$, $v \in N_G^5(u)$ if $d_G(w)=4$, $\alpha(\Gamma_G(w))=2$ for
$w \in \{u,v\}$ and $u$ and $v$ are connected as shown in Figure
\ref{fig:5n}; otherwise $N_G^5(u)=\emptyset$.

For $j=7$, $v \in N_G^7(u)$ if $d_G(w)=4$, $\alpha(\Gamma_G(w))=2$ for
$w \in \{u,v\}$ and $u$ and $v$ are connected as shown in Figure
\ref{fig:7n}; otherwise $N_G^7(u)=\emptyset$.

Note that for $j\in\{1,2,3,5,7\}$, $v\in N_G^j(u)$ if and only if
$u\in N_G^j(v)$; but this does not hold for $j=4$. 
We have the following lemma.

\begin{lemma} \label{N5}
For $u \in V(G)$ such that $d_G(u)=4$ and $\alpha(\Gamma_G(u))=2$,
we have $|N_G^1(u) \cup N_G^2(u) \cup N_G^3(u) \cup N_G^4(u) \cup
N_G^5(u) \cup N_G^7(u)| \leq 96$.
\end{lemma}
{\bf Proof:}  It is clear that  $|N_G^1(u) \cup N_G^2(u) \cup
N_G^3(u)| \leq 4+8+8 \times 3=36$.  We next estimate $|N_G^4(u)|$.
In Figure \ref{fig:4n}, observe that $w$ is
connected to one vertex of $S_2(u)$  and $w \not \in \Gamma_G(u)$.
For a fixed $u$, there are at most four choices for $w$, at most
three choices for $z$, and at most three choices for $v$. Therefore,
we have $|N_G^4(u)| \leq 4 \times 3 \times 3=36$.

Let us estimate $|N_G^5(u)|$. In Figure
\ref{fig:5n},  for a fixed $u$, we have four choices for $w$ and
two choices for $z$. Fix a $z$. Assume $\Gamma_G(z) \setminus
\{w\}=\{a,b,c\}$. Let $T_1=\{a,b\}$, $T_2=\{b,c\}$, and
$T_3=\{a,c\}$. We have the following claim.

\noindent {\bf Claim} There are at most three $v \in N_G^5(u)$ such
that for each $v$  we have $ \Gamma_G(z) \cap \Gamma_G(v)=T_i$ for
some $1 \leq i\leq 3$ as shown in Figure \ref{fig:5n}.

\begin{figure}[htbp]
  \centering{
\psfig{figure=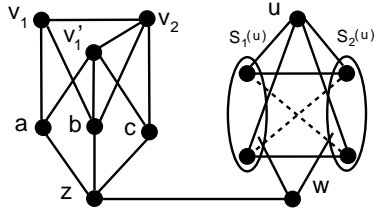, width=0.35\textwidth}}
 \caption{The picture for the claim.}\label{N51}
\end{figure}

\noindent {\bf Proof of the claim:}  For each $1 \leq i \leq 3$,
there are at most three $v \in N_G^5(u)$ such that $ \Gamma_G(z)
\cap \Gamma_G(v)=T_i$ as shown in Figure \ref{fig:5n} since each
vertex in $T_i$ has at most three neighbors other than $z$.  If
the claim is false, then there is  $1 \leq i \not =j \leq 3$ such
that $\Gamma_G(z) \cap \Gamma_G(v_{i})=T_i$ and $\Gamma_G(z) \cap
\Gamma_G(v_{i}')=T_i$  for some $v_{i}, v_{i}' \in N_G^5(u)$,  and
$\Gamma_G(z) \cap \Gamma_G(v_j)=T_j$ for some $v_j \in N_G^5(u)$,
where $v_i,v_i', v_j$ are distinct. Without loss of generality, we
assume $\Gamma_G(z)  \cap \Gamma_G(v_1)=\Gamma_G(z)  \cap
\Gamma_G(v_1')=T_1$ for $v_1, v_1' \in N_G^5(u)$, and  $\Gamma_G(z)
\cap \Gamma_G(v_2)=T_2$ for some $v_2 \in N_G^5(u)$,  see Figure
\ref{N51}. Observe that $\Gamma_G(b)=\{v_1,v_1',v_2,z\}$. Since
$\Gamma_G(z)  \cap \Gamma_G(v_1)=T_1$ as shown in Figure
\ref{fig:5n}, $a$ and one of $b$'s neighbors form $S_i(v_1)$ for
some $i \in \{1,2\}$; we assume it is $S_1(v_1)$. Note
$\{z,v_1,v_1'\} \subset \Gamma_G(a)$. Thus $S_1(v_1)=\{a,v_2\}$ and
$v_2 \in \Gamma_G(v_1)$. Similarly, we can show
$S_1(v_1')=\{a,v_2\}$ and $v_2 \in \Gamma_G(v_1')$.  Now, observe
that $\Gamma_G(v_2)=\{v_1,v_1',b,c\}$. Since $\Gamma_G(z) \cap
\Gamma_G(v_2)=T_2$ as shown in Figure \ref{fig:5n}, $b$ and one of
neighbors of $v_2$ form $S_i(v_2)$ for some $i \in \{1,2\}$; we
assume $i=1$. Because $\{v_1,v_1'\} \subset \Gamma_G(b)$, then
$S_1(v_2)=\{b,c\}$. However, $b$ and $c$ are not is in the same
independent set in the  definition of  $N_G^5(u)$, see Figure
\ref{fig:5n}. This is a contradiction and this case can not
happen.  The claim follows.

Therefore, $|N_G^5(u)| \leq 4 \times 2 \times 3=24.$

In Figure \ref{fig:7n}, for a fixed $u$, we
have two choices for the edge $e$, one choice for $w$, two choices
for $z$, and three choices for the edge $f$.  Fix a $z$. By
considering the degrees of the endpoints of $f$, there is at most
one $f$  and  at most one $v \in N_G^7(u)$ such that $|\Gamma_G(f)
\cap \Gamma_G(v)|=4 $ as shown in Figure \ref{fig:7n}. Therefore,
we have $|N_G^7(u)| \leq 2 \times 2 \times 1=4$.

Last, we estimate $|N_G^5(u) \cup N_G^7(u)|$. If there is some $v
\in N_G^7(u)$, then we observe that there are at most five $z$'s
(see  Figure \ref{fig:5n}). We get the
number of $v \in N_G^5(u)$ is at most $5 \times 3=15$. In this case,
we have
$$
|N_G^5(u) \cup N_G^7(u)| \leq 4+15 < 24.
$$
If $N_G^7(u)=\emptyset$, then also we have
$$
|N_G^5(u) \cup N_G^7(u) \leq 24.$$
 Therefore $$|N_G^1(u) \cup N_G^2(u) \cup N_G^3(u) \cup N_G^4(u)
\cup N_G^5(u) \cup N_G^7(u)| \leq  36+36+24=96.$$ \hfill $\square$

 Based on the graph $G$, we define an auxiliary   graph $G^{\ast}$ on
vertex set $V(G)$. The edge set is defined as follows:  $uv \in E(G^{\ast})$ precisely if
 either $u \in N_G^1(v) \cup N_G^2(v) \cup N_G^3(v) \cup N_G^4(v) \cup N_G^5(v)
\cup N_G^7(v) $, or $v \in N_G^4(u)$.
 We have the following lemma.

\begin{lemma} \label{121color}
The graph $G^{\ast}$ is 133-colorable.
\end{lemma}

\noindent
{\bf Proof:}  Let $\sigma$ be an increasing order of $V(G^{\ast})$
satisfying the following conditions.
\begin{description}
\item [1:] For $u$ and $v$ such that  $d_G(u)=3$ and $d_G(v)=4$, we have $\sigma(u) <
\sigma(v)$.
\item [2:]  For $u$ and $v$ such that $d_G(u)=d_G(v)=4$, $\alpha(\Gamma_G(u)) \geq
3$, and $\alpha(\Gamma_G(v))=2$, we have $\sigma(u) < \sigma(v)$.
\end{description}
We will color $V(G^{\ast})$ according to the order $\sigma$. For
each $v$, we have the following estimate on the number of colors
forbidden to use for $v$.
\begin{description}
\item [1:] For $v$ such that $d_G(v)=3$, the number of colors
forbidden to use for $v$ is at most $|N_G^1(v) \cup N_G^2(v) \cup
N_G^3(v)| \leq 3+9+27=39$.
\item [2:] For $v$ such that $d_G(v)=4$ and $\alpha(\Gamma_G(v)) \geq
3$,  the number of colors forbidden to use for $v$ is at most
$|N_G^1(v) \cup N_G^2(v) \cup N_G^3(v)| \leq 3+9+27=39$.
\item [3:]  For $v$ such that $d_G(v)=4$ and
$\alpha(\Gamma_G(v))=2$, the number of colors forbidden to use for
$v$ is at most  $|N_G^1(v) \cup N_G^2(v) \cup N_G^3(v) \cup N_G^4(v)
\cup N_G^5(v) \cup N_G^7(v)|+|N_G^4(v)| \leq 96+36=132$ by Lemma
\ref{N5}.
\end{description}
Therefore, the greedy algorithm shows $G^{\ast}$ is 133-colorable.
\hfill $\square$

Let $X$ be a color class of $G^{\ast}$. We define a new graph
$G(X)$ by the following process.

\begin{enumerate}
\item For each $x \in X$,  if $|S(x)|=2$ or $|S(x)|=3$, then we contract $S(x)$ as a single vertex,
delete the vertices in $\Gamma_G(v) \setminus S(v)$, and keep
label 1 on the new vertex;  if $|S(x)|=4$, i.e., $S(x)=S_1(x)
\cup S_2(x)$, then we contract $S_1(x)$ and $S_2(x)$ as single
vertices and keep their labels. After that, we  delete $X$. Let $H$
be the resulting graph.
\item Note that $\Gamma_H(x) \cap \Gamma_H(y) = \emptyset$  and   there is no edge from
 $\Gamma_H(x)$ to  $\Gamma_H(y)$ for any $x,y \in X$ as $X$ is a color class.
\item  We identify all vertices with label $i$ as a single vertex $w_i$ for $i \in \{1,2\}$. Let $G(X)$ be the resulted graph.
\end{enumerate}
We have the following lemma on the chromatic number of $G(X)$.
\begin{lemma} \label{4colorable}
The graph $G(X)$ is $4$-colorable for each color class.
\end{lemma}
We postpone the proof of this lemma until the end of this section
and prove Lemma  \ref{f4} first.

\noindent
{\bf Proof of Lemma \ref{f4}:} By Lemma \ref{121color},  there is a
proper 133-coloring of $G^{\ast}$. We assume
$V(G^{\ast})=V(G)=\cup_{i=1}^{133} X_i $, where $X_i$ is the $i$-th
color class.

For each $i \in \{1,\ldots,  133\}$,  Lemma \ref{4colorable} shows
$G(X_i)$ is $4$-colorable; let $c_i\colon V(G(X_i)) \to T_i$ be a
proper $4$-coloring of the graph $G(X_i)$. Here $T_1,T_2,\ldots,
T_{133}$ are pairwise disjoint; each of them consists of $4$ colors.
For $i \in \{1,\ldots,  133\}$, the $4$-coloring $c_i$ can be viewed
as a $4$-coloring of $G \setminus X_i$ since each vertex with label
$j$ receives the color $c_i(w_i)$ for $j=1,2$ and  each removed
vertex has at most three neighbors in $G \setminus X_i$.

Now we  reuse the notation $c_i$ to denote this $4$-coloring of $G
\setminus X_i$. For each $v \in X_i$, we have $|\cup_{u \in
\Gamma_G(v)}c_i(u)| \leq  2$. We can assign two unused colors,
denoted by the set $Y(v)$, to $v$. We  define $f_i: V(G) \rightarrow
 \mathcal P(T_i)$ (the power set of $T_i$) satisfying
 $$
 f_i(v)= \left\{
 \begin{array}{ll}
 \{c_i(v)\}  \   & \textrm{if} \  v \in V(G) \setminus X_i, \\
 Y(v) &  \textrm{if} \  v \in X_i.\\
 \end{array}
 \right.
 $$

 Observe that   each vertex in $X_i$ receives two colors from $f_i$ and
  every other vertex receives one color.  Let $\sigma: V(G)
 \rightarrow \mathcal P(\cup_{i=1}^{k} T_i)$ be
a mapping  such that  $\sigma(v)=\cup_{i=1}^{m} f_i(v)$. It is easy
to verify
 $\sigma$ is a $134$-fold coloring of $G$ such that each color is drawn
 from a palette of $532$ colors; namely we have
 $$\chi_f(G) \leq
 \frac{532}{134} =4-\frac{2}{67}.$$
The proof of Lemma \ref{f4} is finished. \hfill $\square$

Before we prove Lemma \ref{4colorable}, we need the following definitions.

A {\it block} of a graph is a maximal $2$-connected induced subgraph.
A {\it Gallai tree} is a connected graph in which all blocks are
either complete graphs or odd cycles. A {\it Gallai forest} is a
graph all of whose components are Gallai trees. A $k$-{\it Gallai
tree (forest)} is a Gallai tree (forest) such that the degree of all
vertices are at most $k-1$. A $k$-critical graph is a graph $G$
whose chromatic number is $k$ and deleting any vertex can decrease
the chromatic number. Gallai showed the following Lemma.

\begin{lemma}\label{gallai} \cite{gallai}
If $G$ is a $k$-critical graph, then the subgraph of $G$ induced on
the vertices of degree $k-1$ is a k-Gallai forest.
\end{lemma}

\noindent
 {\bf Proof of Lemma \ref{4colorable}:} We use proof by
 contradiction.  Suppose that $G(X)$ is not $4$-colorable.
 The only possible vertices in $G(X)$ with degree greater than $4$ are the
 vertices $w_1$ and $w_2$, which are
 obtained by contracting the vertices with label  1 and 2 in the intermediate graph $H$.
The simple greedy algorithm shows that $G(X)$ is always
$5$-colorable. Let
 $G'(X)$ be a $5$-critical subgraph of $G(X)$.  Applying
 Lemma \ref{gallai} to $G'(X)$, the subgraph of $G'(X)$ induced on the
 vertices of degree $4$ is a $5$-Gallai forest $F$. The vertex set of $F$
may contain  $w_1$ or $w_2$. Delete $w_1$ and $w_1$ from $F$ if $F$
contains one of  them. Let $F'$ be the resulting Gallai forest. (Any
induced subgraph of a Gallai forest is still a Gallai forest.) The  Gallai forest $F'$ is not empty. Let $T$ be a
connected component of $F'$ and $B$ be a leaf block of $T$. The
block $B$ is either a clique or an odd cycle from the definition of a
Gallai tree.

Let $v$ be a vertex in $B$. As $v$ has at most two neighbors ($w_1$
and $w_2$) outside $F'$ in $G(X)$,   we have $d_{F'}(v)\geq 2$. If
$v$ is not in other blocks of $F'$, then we have $d_B(v)\geq 2$.
It follows that $|B|\geq 3$. Since $B$ is a subgraph of $G$ and
$G$ is $K_4$-free, the block $B$ is  an
odd cycle.

Let $v_1v_2$ be an edge in $B$ such that $v_1$ and $v_2$ are not in
other blocks. The degree requirement implies $v_iw_j$ are edges in
$G(X)$ for all $i, j\in\{1,2\}$. For $i=1,2$, there are vertices
$x_i, y_i\in X$ satisfying $S(x_i)\cap \Gamma_G(v_i)\not =\emptyset$
and $S(y_i)\cap \Gamma_G(v_i)\not =\emptyset$; moreover either
$|S(x_i)| =4$ or $|S(y_i)|=4$ since one of its neighborhood has
label 2. Without loss of generality, we assume $|S(x_i)|=4$ for $i
\in \{1,2\}$. If $x_i\not= y_i$, then $y_i \in N_G^4(x_i)$, i.e.,
$y_i \in \Gamma_{G^{\ast}}(x_i)$; this contradicts $X$ being a
color class. Thus we have $x_i=y_i$
 and $|S(x_i)|=4$ for $i \in \{1,2\}$. For
$\{i,j\}=\{1,2\}$, if $x_i\not= y_j$,  then $y_i \in N_G^5(x_i)$,
 i.e., $y_i \in \Gamma_{G^{\ast}}(x_i)$; this is a contradiction of $X$ being a color class.
 Thus we
have
$$x_1=x_2 = y_1=y_2.$$
Let $x$ denote this common vertex above. Then
$d_G(x)=4$ and $\alpha(\Gamma_G(x))=2$.

Let $v_0$ be  the only vertex in $B$ shared by other blocks. Since
$B-v_0$ is connected, the argument above shows there is a common $x$
for all edges in $B-v_0$. If $\Gamma_{G(X)}(v_0)\cap
\{w_1,w_2\}\not=\emptyset$, the there is some vertex $x_0 \in X$
such that $S(x_0) \cup \Gamma_G(v_0) \not =\emptyset$. By the
similar argument, we also have $x_0=x$.

Therefore, $x$ depends only on $B$. In the sense that for any $y\in X$ and any $v\in B$, if $S(y)\cap \Gamma_G(v)\not=\emptyset$,
then $y=x$.

 The block $B$ is an odd cycle as we mentioned above.  Suppose $|B|=2r+1$. Let $v_0,v_1,\ldots,v_{2r}$ be the vertices
of $B$ in cyclic order and $v_0$  be the only vertex which may be shared by other block.

Let $x\in X$ be the vertex determined by $B$. Recall
$d_G(x)=4$ and $\alpha(\Gamma_G(x))=2$.
Each vertex in $\Gamma(x)$ can have at most $2$ edges to $B$.  We get
\begin{equation}
  \label{eq:12}
4r\leq |E(B,\Gamma(x))|\leq 8.
\end{equation}
We have $r\leq 2$. The block $B$ is either a $C_5$ or a $K_3$.
We claim both $v_0w_1$ and $v_0w_2$ are non-edges of
$G(X)$.

If $B=C_5$, then inequality (\ref{eq:12}) implies
 that $v_0$ has no neighbor in
$\Gamma(x)$ and the claim holds.
If $B=K_3$, then the claim also holds; otherwise $B\cup
\{\underline{S_1(x)}, \underline{S_2(x)}\}$ forms a $K_5^-$ in
$G/S_1(x)/S_2(x)$, which is a contradiction to Lemma \ref{D41}.

 Let $u_1$ and $u_2$ be the two neighbors
of $v_0$ in other blocks of $F'$. If $u_1$ and $u_2$ are in the same block,
then this block is an odd cycle; otherwise, $v_0u_1$ and $v_0u_2$ are in two different blocks.

The union of non-leaf blocks of $T$ is a Gallai-tree, denoted by
$T'$. The argument above shows every leaf block of $T'$ must be an
odd cycle. Let $C$ be such a leaf block of $T'$.  Now $C$ is an odd
cycle, and $C$ is connected to $|C|-1$ leaf blocks of $T$. Let $B$
and $B'$ be two leaf blocks of $T$ such that $B\cap C$ is adjacent
to $B'\cap C$. Without loss of generality, we may assume $B$ is the
one we considered before.  By the same argument, $B'$ is an odd
cycle of size $2r'+1$ with $r'\in \{1,2\}$.  Let $v'_0,v'_1,\ldots,
v'_{2r'}$ be the vertices of $B'$ and $v'_0$ be the only vertex in
$B'\cap C$. For $i$ in $\{1,2,\ldots, 2r'\}$ and $j$ in $\{1,2\}$,
$v'_iw_j$ are edges in $G'(X)$. Similarly,  there exists a vertex
$x'\in X$ with $d_G(x')=4$ and $\alpha(\Gamma_G(x'))=2$ such that
$|E(v_i, S_1(x'))|\geq 1$ and $|E(v_i, S_2(x'))|\geq 1$.  We must
have $x=x'$; otherwise $x'\in N_G^7(x)$, i.e., $x' \in
\Gamma_{G^{\ast}}(x)$, and this contradicts the fact that $X$ is a
color class in $D$. Now we have $|E(\Gamma(x), B)|\geq 4r$ and
$E(\Gamma(x), B')|\geq 4r'$. By counting  the degrees of vertices in
$\Gamma(x)$ in $G$, we have
$$4r+4r'+4+4\leq 16.$$
We get  $r=r'=1$. Both $B$ and $B'$ are $K_3$'s. In this case,
$G/S_1(x)/S_2(x)$ contains the graph $G_0$, see figure
\ref{fig:H13}. This contradicts Lemma \ref{D41}.

We can find the desired contradiction, so the lemma
follows.  \hfill $\square$

\section{Proof of  Lemma \ref{D41}}

In this section, we will prove  Lemma \ref{D41}. We first review a Lemma from \cite{lup}.

\begin{lemma}\label{cut2}
Let $G$ be a graph.  Suppose that $G_1$ and $G_2$ are two subgraphs
such that $G_1 \cup G_2=G$ and $V(G_1) \cap V(G_2) =\{u,v\}$.
\begin{enumerate}
\item If $uv$ is an edge of $G$, then we have
$$\chi_f(G) = \max\{\chi_f(G_1), \chi_f(G_2)\}.$$
\item If $uv$ is not an edge of $G$, then we have
$$\chi_f(G) \leq \max\{\chi_f(G_1), \chi_f(G_2+uv), \chi_f(G_2/uv)\},$$
\end{enumerate}
where $G_2+uv$ is the graph obtained from $G_2$ by adding edge $uv$ and $G_2/uv$ is the graph obtained from $G_2$ by contracting  $\{u, v\}$.
\end{lemma}

\noindent {\bf Proof of Lemma \ref{D41}:} Recall that $G$ is a
connected $K_4$-free graph with minimum number of vertices such that
$G \not = C_8^2$ and $ \chi_f(G) > 4-\frac{2}{67} $. Note that $G$
is 2-connected.  We will prove it by contradiction.

\begin{figure}[htbp]
 \centerline{\psfig{figure=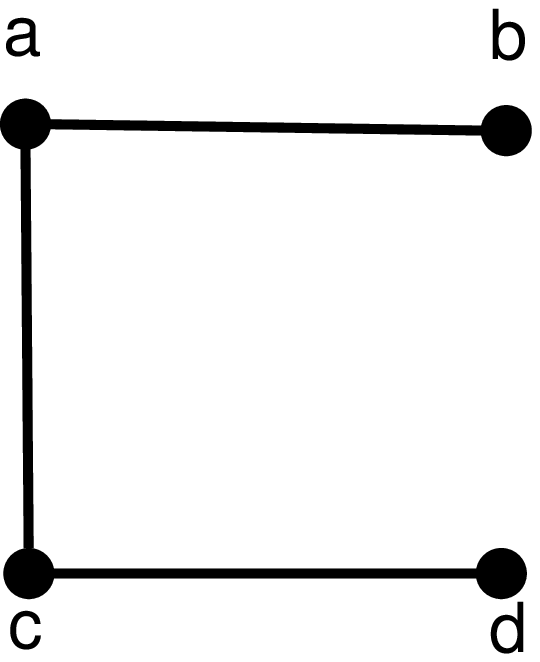,width=0.18\textwidth}  \hfil
 \psfig{figure=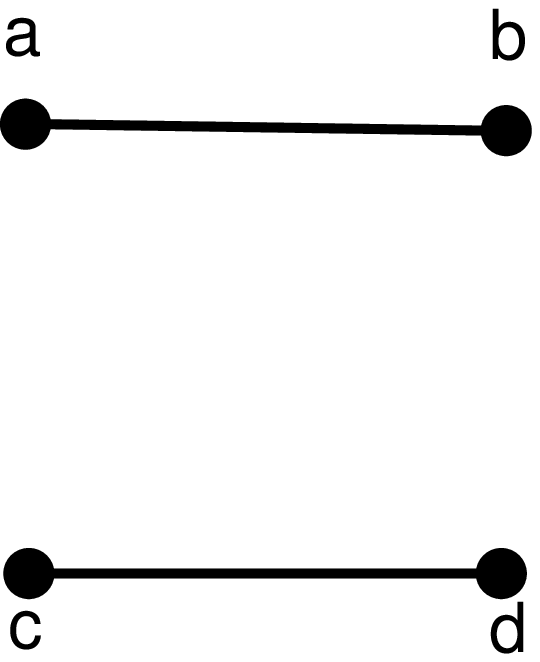, width=0.18\textwidth} \hfil   \psfig{figure=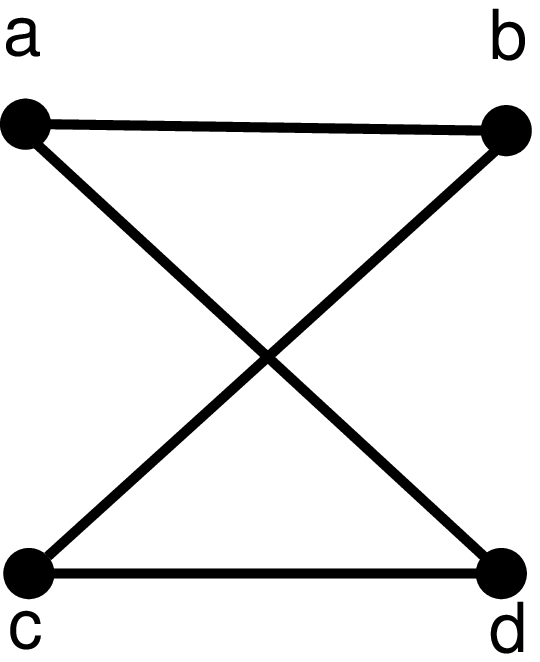, width=0.18\textwidth}
}
 \centerline{$\Gamma_G(x)=P_4$   \hfil $\Gamma_G(x)=2~e$ \hfil $\Gamma_G(x)=C_4$}
\caption{Three possible cases of $\Gamma_G(x)$. } \label{Delta=4}
\end{figure}

Suppose Lemma \ref{D41} fails for some vertex $x$ in $G$. Observe
$\Gamma_G(x)$ is one of the graphs in  Figure \ref{Delta=4}. Here we
assume $\Gamma_G(x)=\{a,b,c,d\}$. Through the proof of the lemma, let
$S_1$ and $S_2$ be two vertex-disjoint independent sets in
$\Gamma_G(x)$, $H$ be a triangle in $V(G) \setminus ( \{x\} \cup
\Gamma_G(x))$, then say $(S_1,S_2,H)$ is a bad triple if
$\{\underline {S_1}, \underline {S_2}, H\}$ contains a $K_5^-$ in
$G/S_1/S_2$.

\begin{figure}[htbp]
 \centerline{\psfig{figure=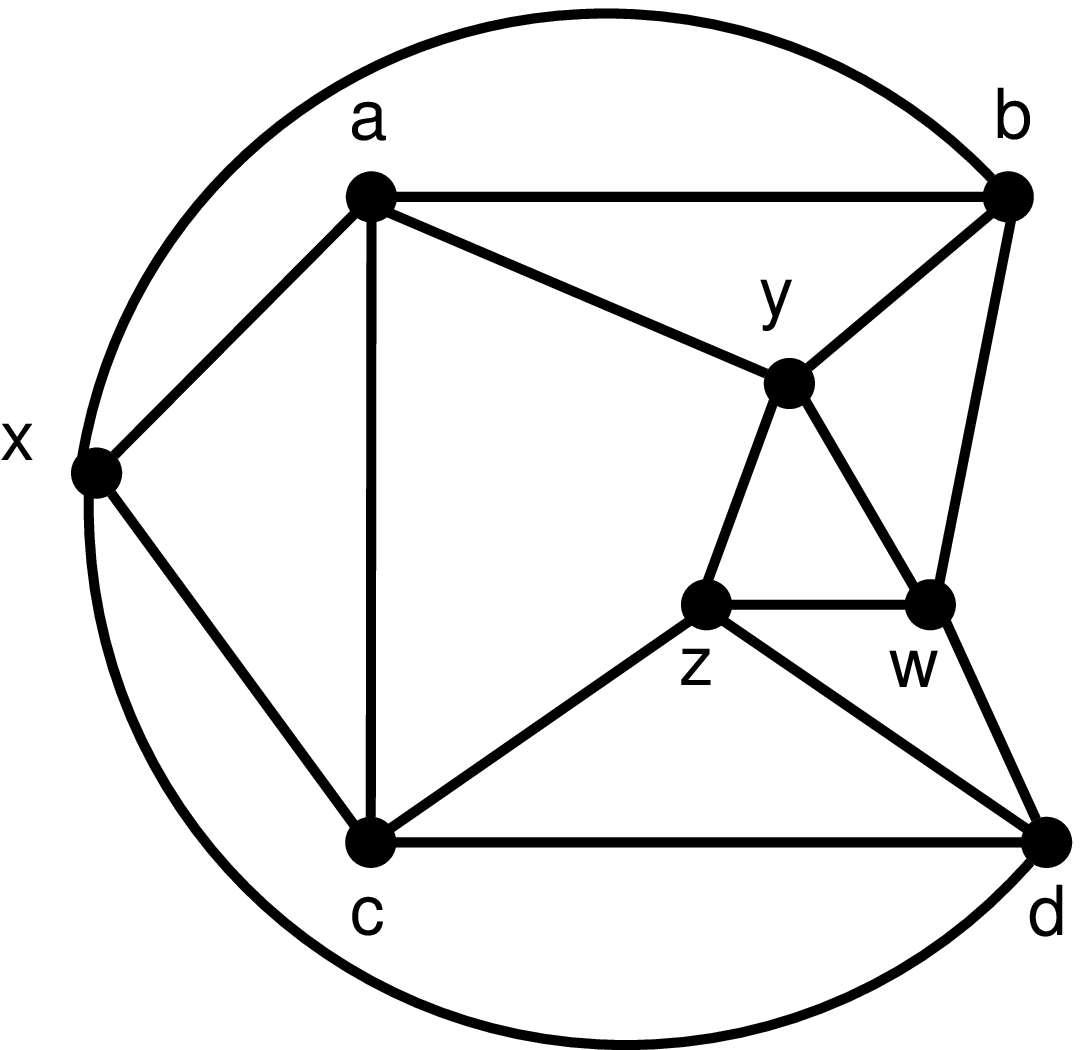, width=0.25\textwidth} \hfil
 \psfig{figure=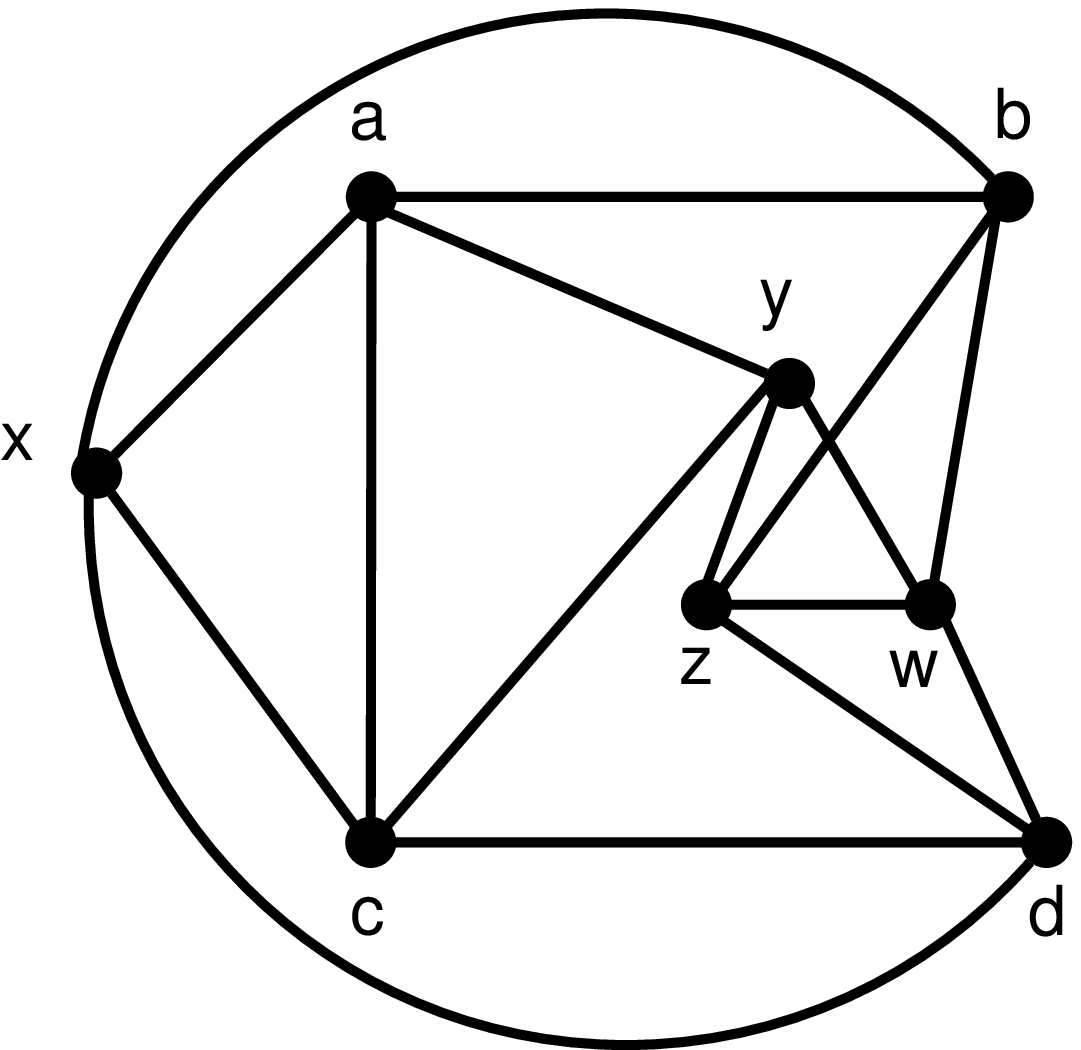, width=0.25\textwidth} \hfil
 \psfig{figure=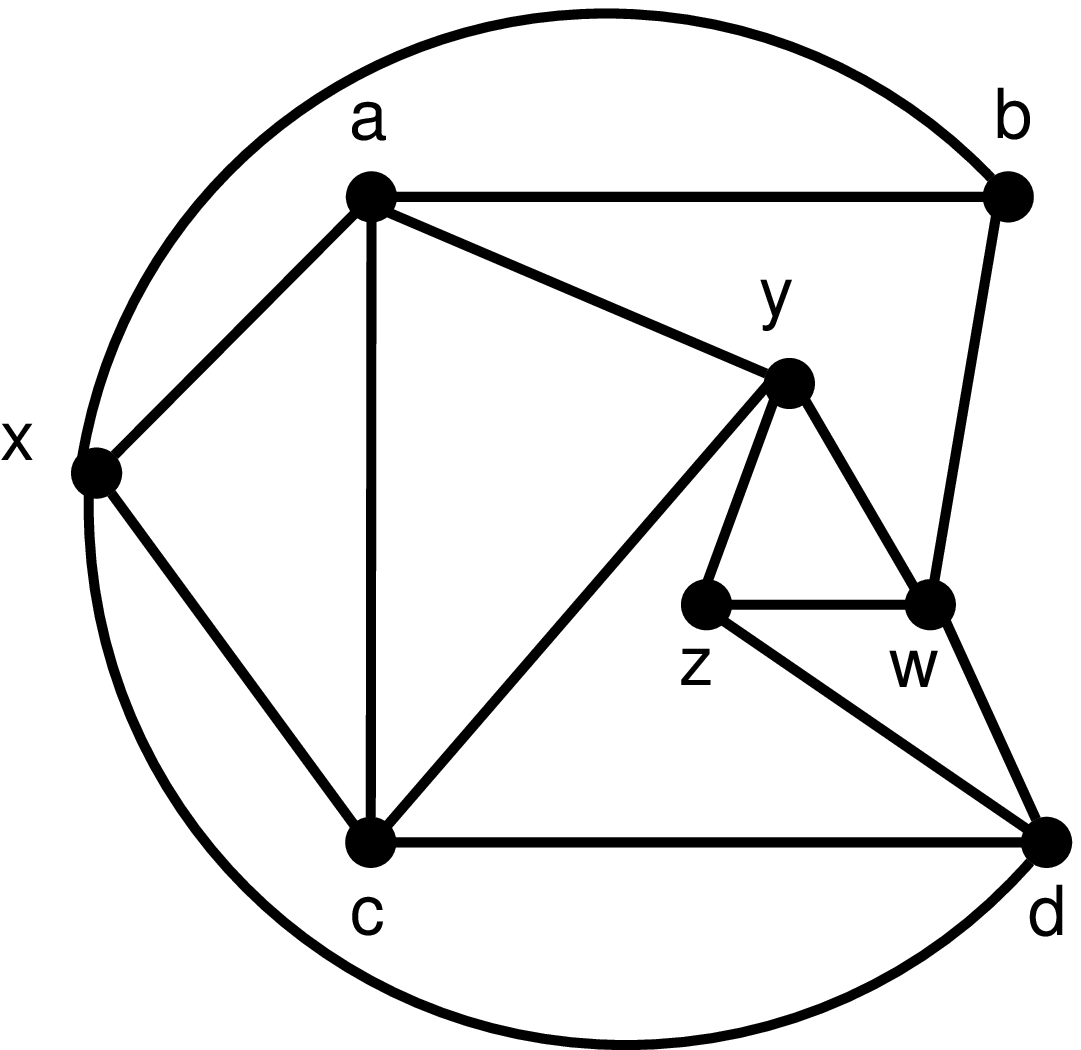, width=0.25\textwidth} }
\centerline{ \hspace*{0.05\textwidth} $H_{1}$ \hspace*{0.15\textwidth} \hfil $H_{2}$ \hspace*{0.15\textwidth}  \hfil $H_{3}$ }

\vspace*{5mm}
 \centerline{\psfig{figure=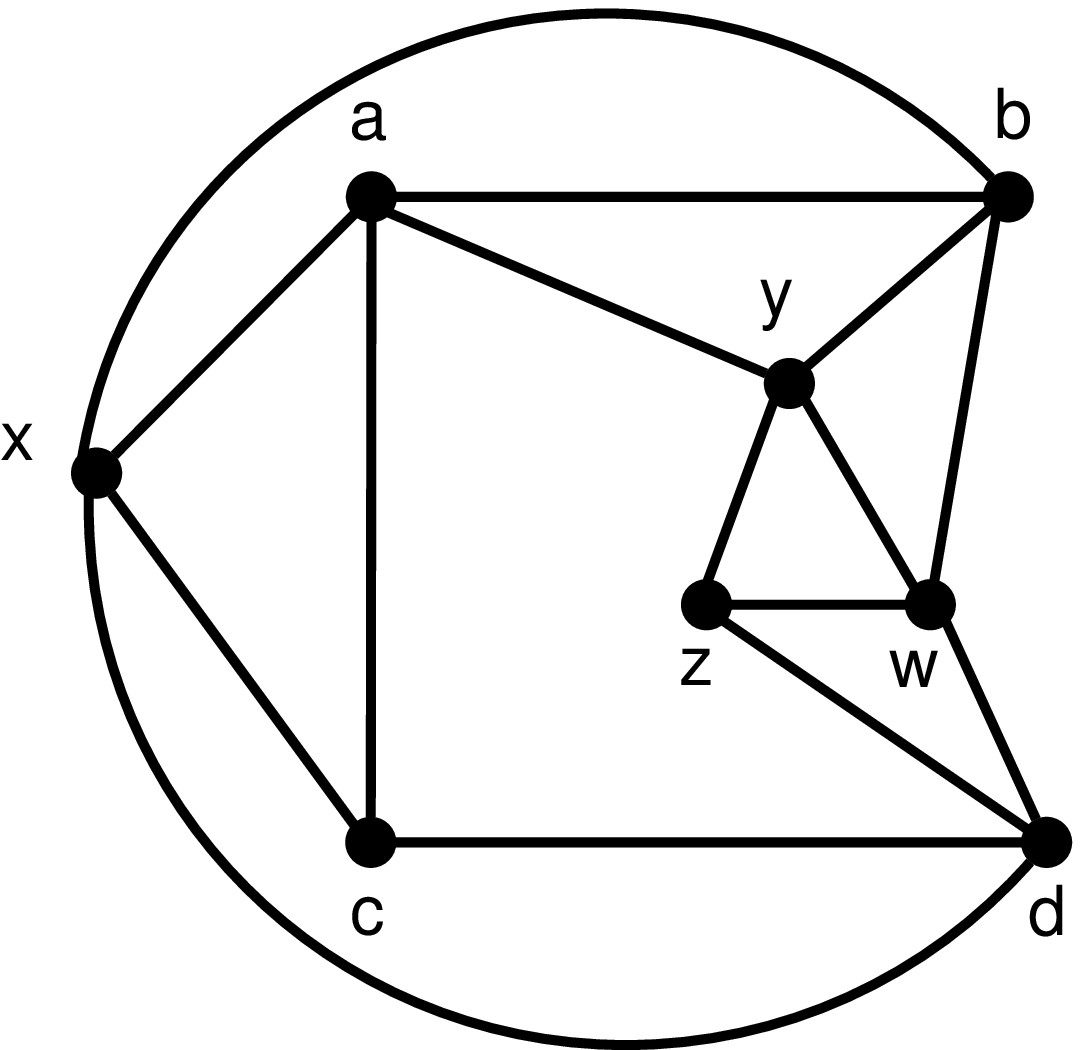,width=0.25\textwidth} \hfil
\psfig{figure=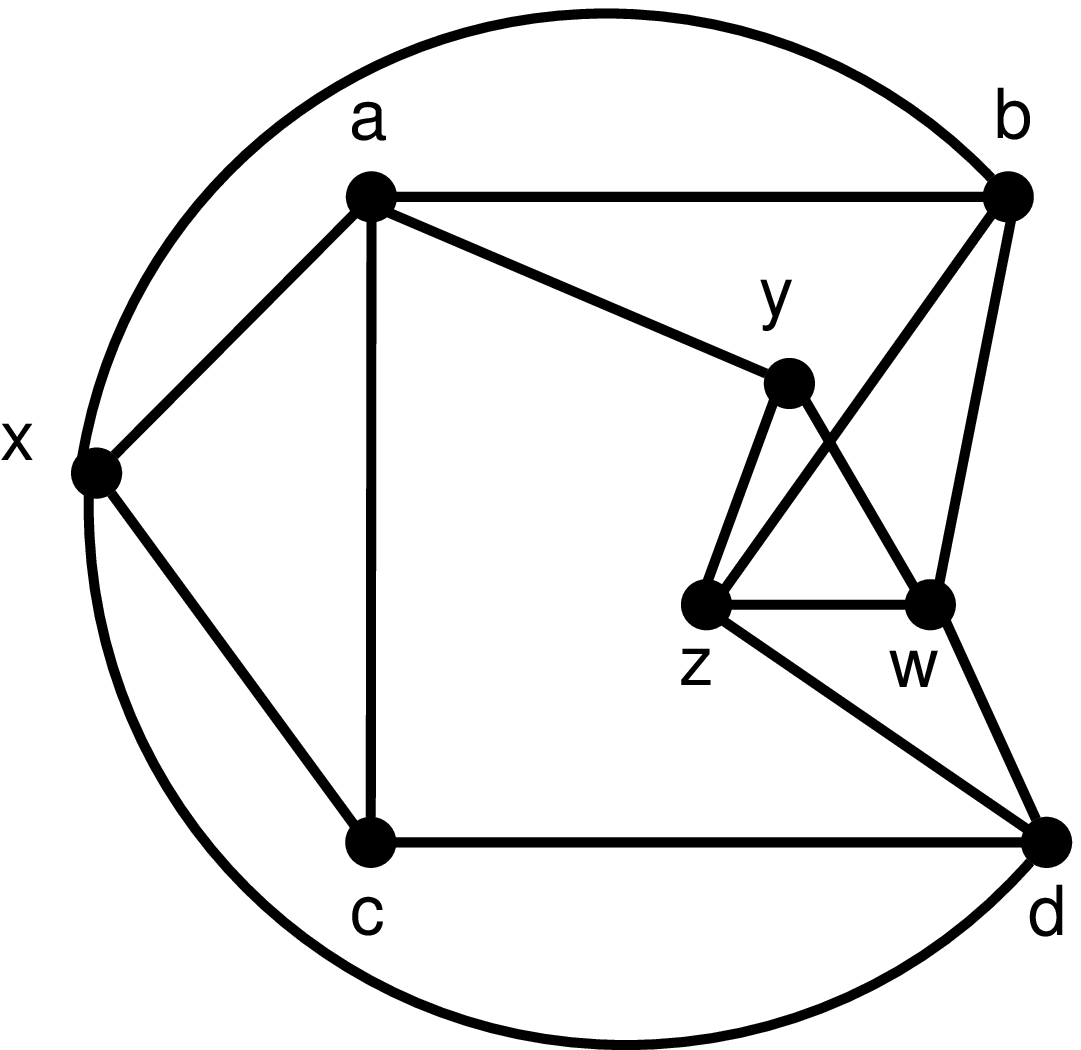, width=0.25\textwidth}
 \hfil
\psfig{figure=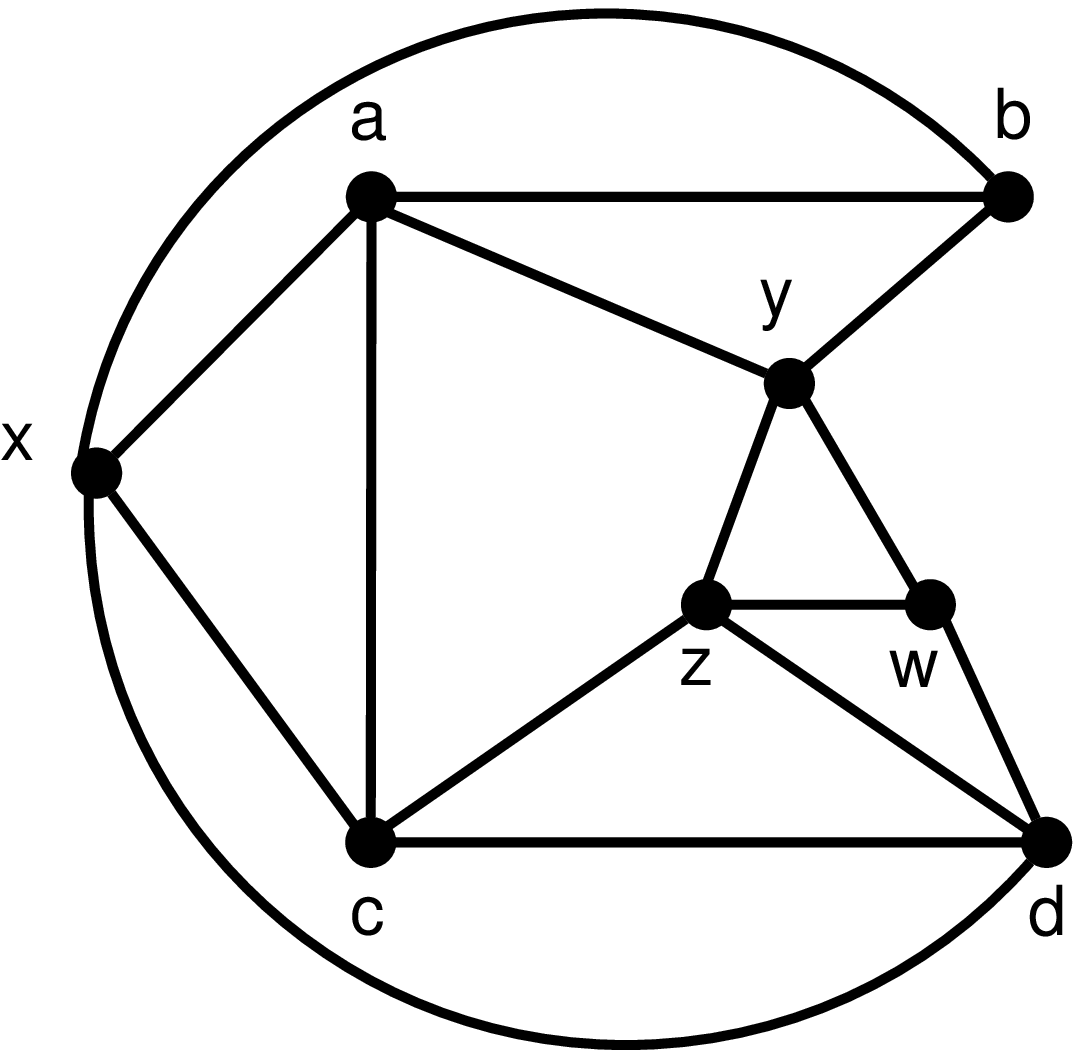,
 width=0.25\textwidth}}
\centerline{  \hspace*{0.05\textwidth}$H_{4}$  \hfil  \hspace*{0.15\textwidth} $H_{5}$  \hfil  \hspace*{0.15\textwidth} $H_{6}$ }
 \caption{If $\Gamma_G(x)=P_4$, then there are six possible induced subgraphs.}
 \label{fig:P4}
\end{figure}

If $\Gamma_G(x)=P_4$, then $\{a,d\}$ and $\{b,c\}$ is the only pair
of disjoint non-edges. There is a triangle  $H$ with
$V(H)=\{y,z,w\}$ such that $(\{a,d\}, \{b,c\}, H)$ is a bad triple.
Note that $|E(\{a,b,c,d\},\{y,z,w\})|=5$ or $6$.
 By an exhaustive search, the induced subgraph of $G$ on $\{x, a,b,c,d,y,z,w\}$
is  one of the following six graphs (see Figure \ref{fig:P4}).

\begin{figure}[htbp]
 \centerline{ \psfig{figure=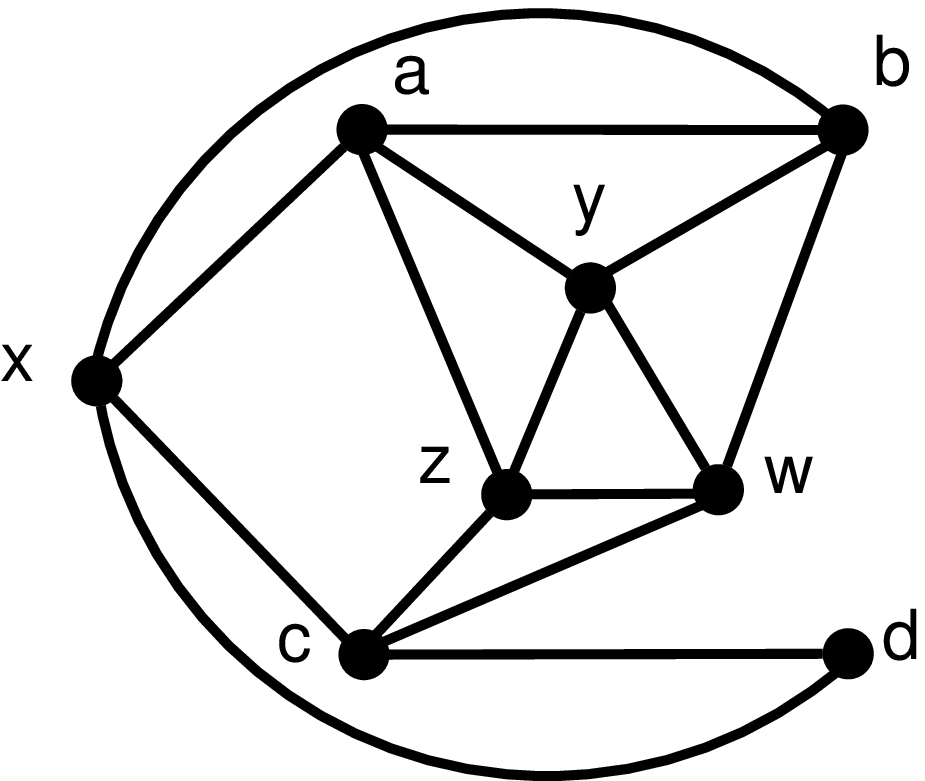, width=0.25\textwidth} \hfil  \hfil \psfig{figure=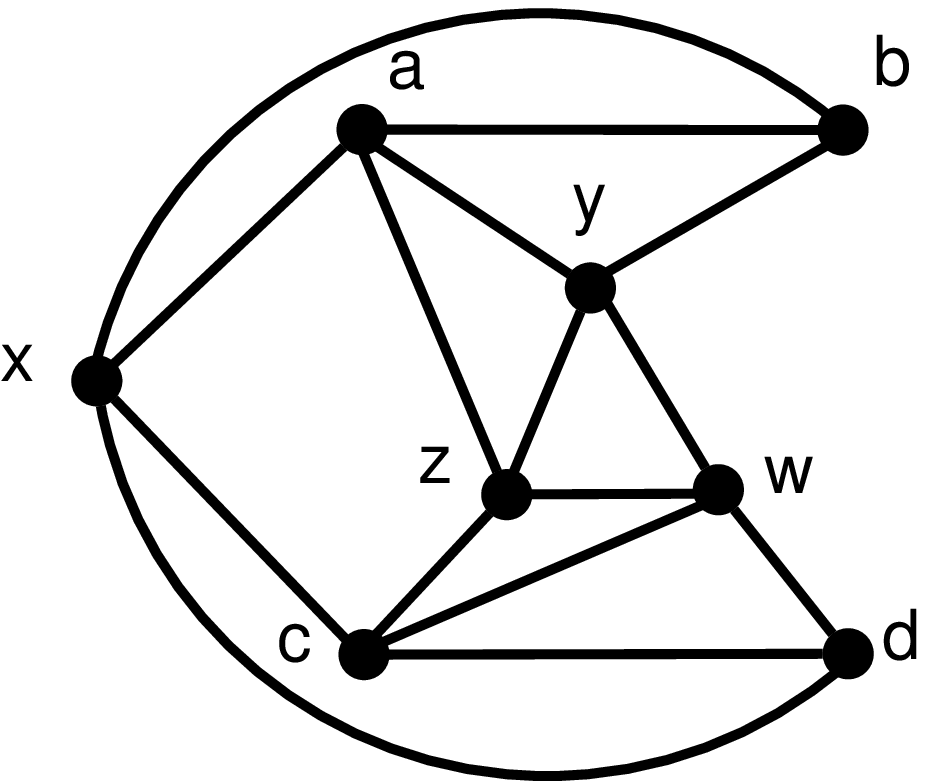,
 width=0.25\textwidth}  \hfil \psfig{figure=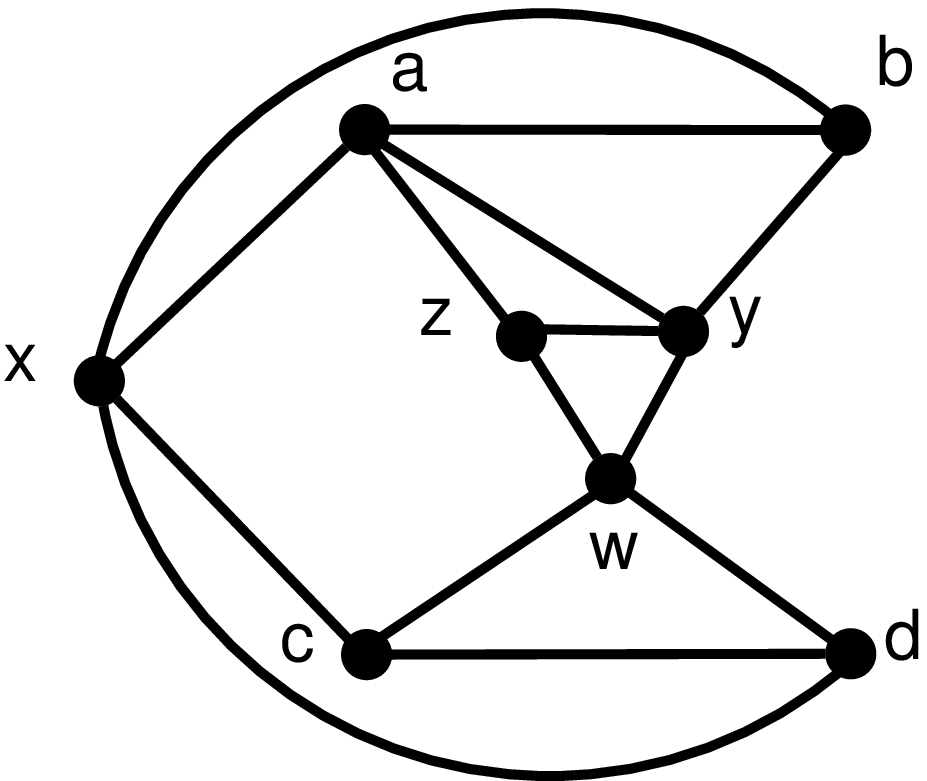,width=0.25\textwidth} }
 \centerline{ \hspace*{0.05\textwidth} $H_7$  \hspace*{0.15\textwidth} \hfil $H_8$   \hspace*{0.15\textwidth} \hfil $H_9$ }
 \caption{If $\Gamma_G(x)=2\ e$, then there are three possible induced subgraphs.}
 \label{fig:C4}
\end{figure}

If $\Gamma_G(x)=2~e$, then $(\{a,c\}, \{b,d\})$ and $(\{a,d\},
\{b,c\})$ are two pairs of disjoint non-edges. By considering the
degrees of vertices in $\Gamma_G(x)$,  there is only one triangle
$H$ with $V(H)=\{y,z,w\}$ such that $(\{a,c\}, \{b,d\}, H)$ and
$(\{a,d\}, \{b,c\}, H)$ are two bad triples. By an exhaustive
search, the induced subgraph of $G$ on $\{x, a,b,c,d,y,z,w\}$ is one
of the following three graphs (see Figure \ref{fig:C4}).

\begin{figure}[htbp]
  \centering
  \psfig{figure=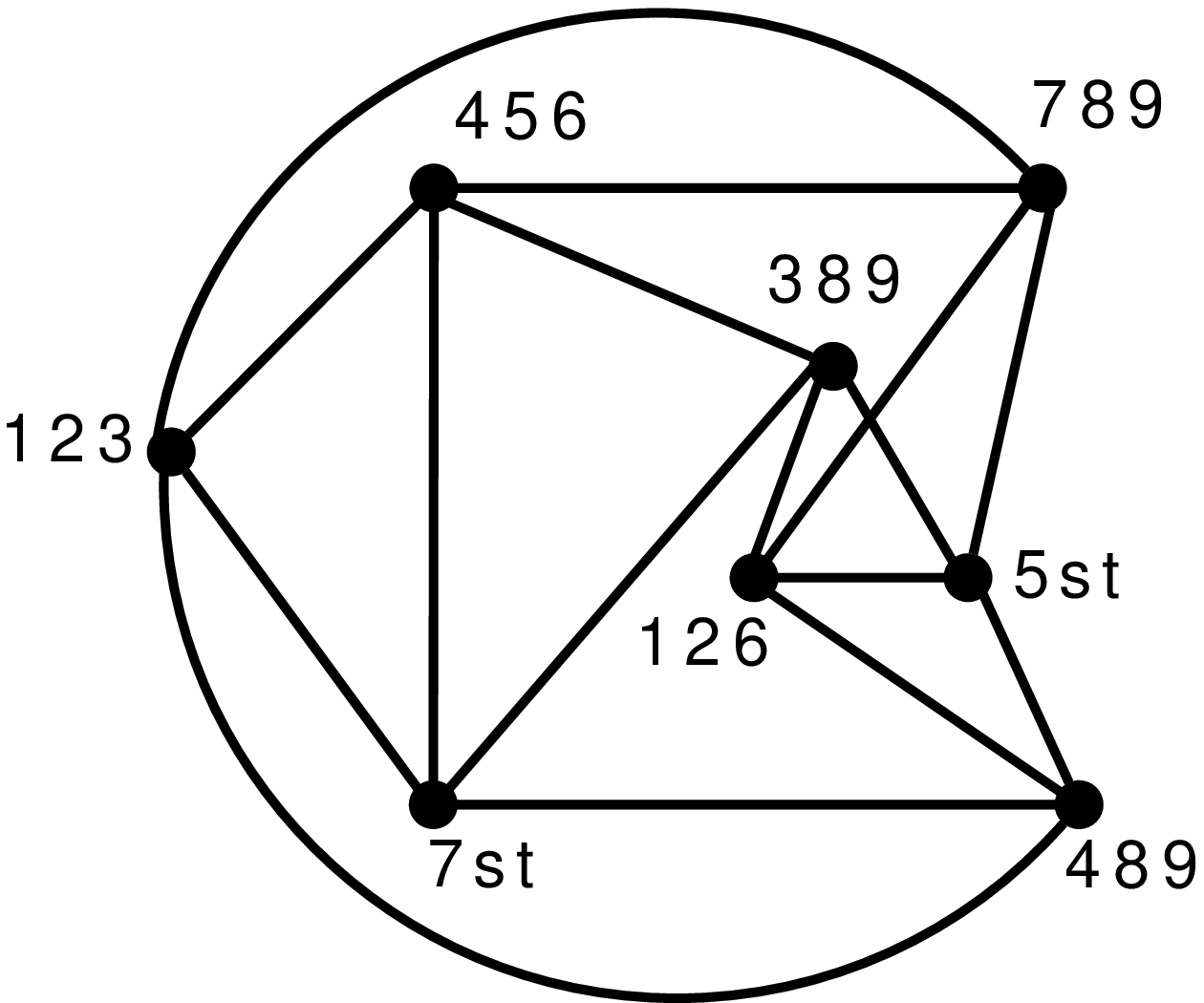, width=0.24\textwidth} \hfil
 \psfig{figure=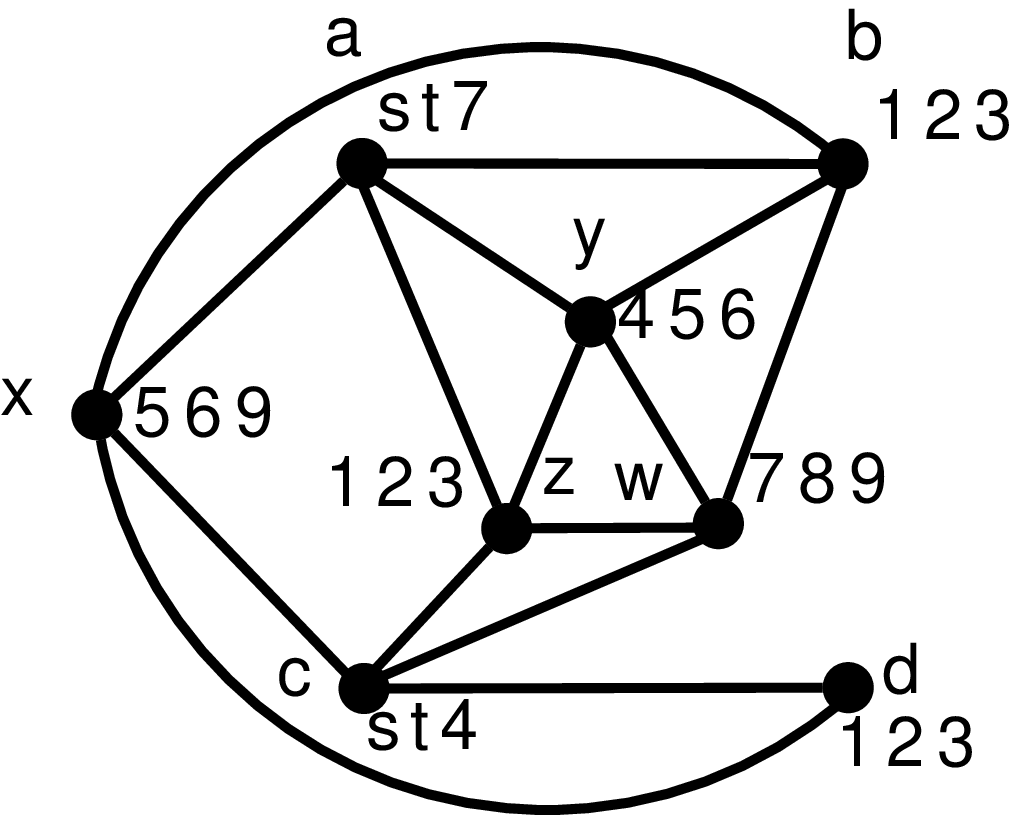, width=0.24\textwidth}
\centerline{$H_2$  \hspace*{1cm} \hfil $H_7$  }
  \caption{The graph $H_2$ and $H_7$ are 11:3-colorable.}
  \label{fig:H2H7}
\end{figure}

It suffices to show that $G$ cannot contain $H_i$ for $1\leq i\leq
9$. Since all vertices in $H_1$ (and $H_2$) have degree 4,  $H_1$
(and $H_2$) is the entire graph $G$. Observe that $H_1$ is
isomorphic to $C_8^2$ and $H_2$ is 11:3-colorable (see Figure
\ref{fig:H2H7}). Contradiction!

In $H_7$, the vertex $d$ is the only vertex with degree less than
$4$. If $H_7$ is not the entire graph $G$, then $d$ is a cut vertex
of $G$. This contradicts the fact that $G$ is $2$-connected. Thus
$G=H_7$. The graph $H_7$ is 11:3-colorable as shown by Figure
\ref{fig:H2H7}. Contradiction!

Now we consider the case $H_3$. Note $H_3+bz$ is the graph $H_2$. We
have $\chi_f(H_3)\leq \chi_f(H_2)\leq 11/3$. The graph $H_3$ must be
a proper induced subgraph of $G$, and the pair $\{b,z\}$ is a vertex
cut of $G$. Let $G'$ be the induced subgraph of $G$ by deleting all
vertices in $H_3$ but $b,z$. We apply Lemma \ref{cut2} to $G+bz$
with $G_1=H_3+bz=H_2$ and $G_2=G'+bz$. We have
$$\chi_f(G+bz)\leq \max\{\chi_f(H_2), \chi_f(G'+bz)\}.$$
Note $\chi_f(H_2)\leq 11/3$ and $11/3<\chi_f(G)\leq \chi_f(G+bz)$.
We have $\chi_f(G)\leq \chi_f(G'+bz)$. Both $b$ and $z$ have at most $2$ neighbors
in $G'+bz$. Thus $G'+bz$ is $K_4$-free; $G'+bz\not=C_8^2$
and has fewer vertices than $G$.
This contradicts to the minimality of $G$.

Note $H_5+cy=H_2$.  The case $H_5$ is similar to the case $H_3$.

\begin{figure}[htbp]
  \centering
  \psfig{figure=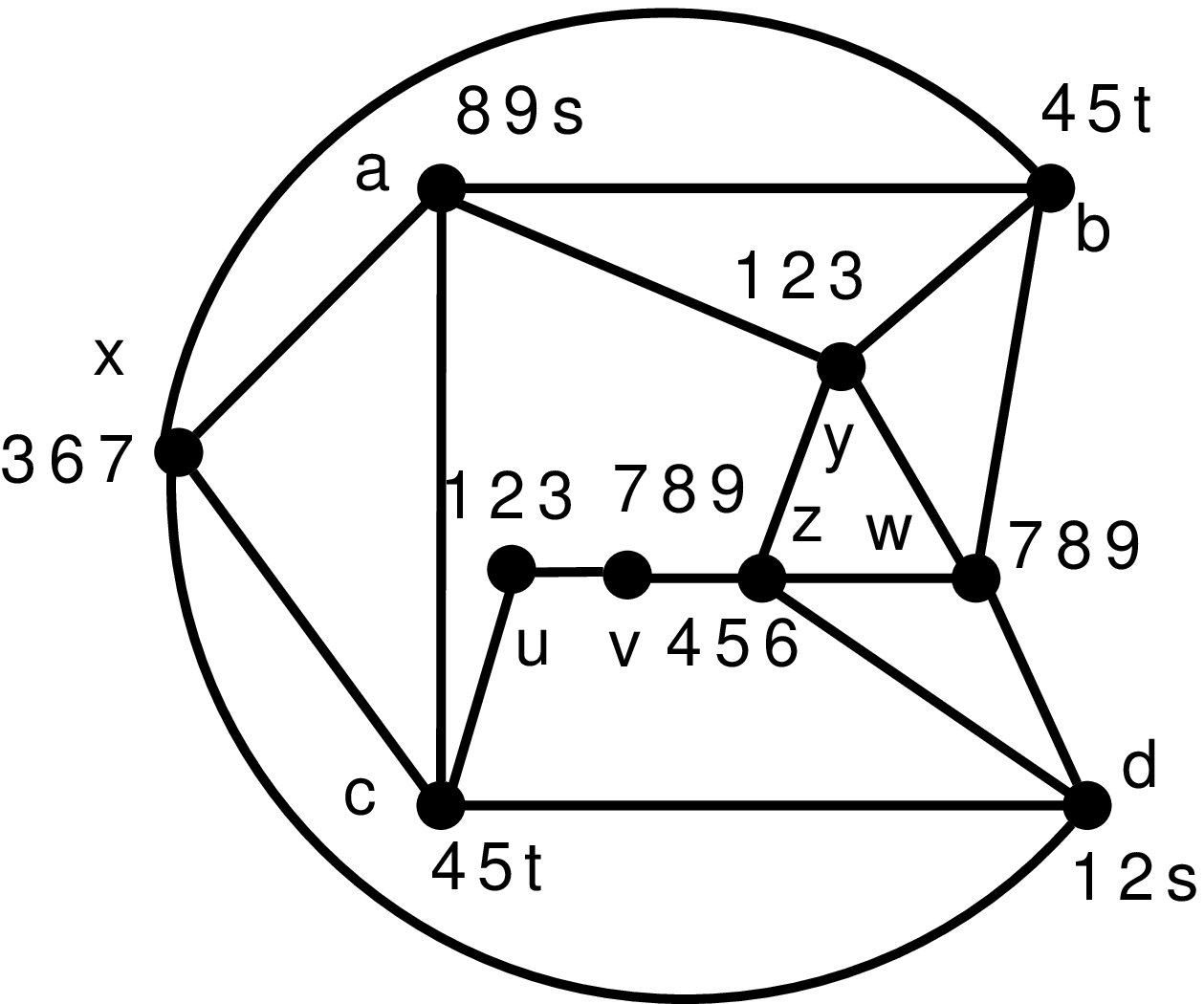, width=0.24\textwidth} \hfil
 \psfig{figure=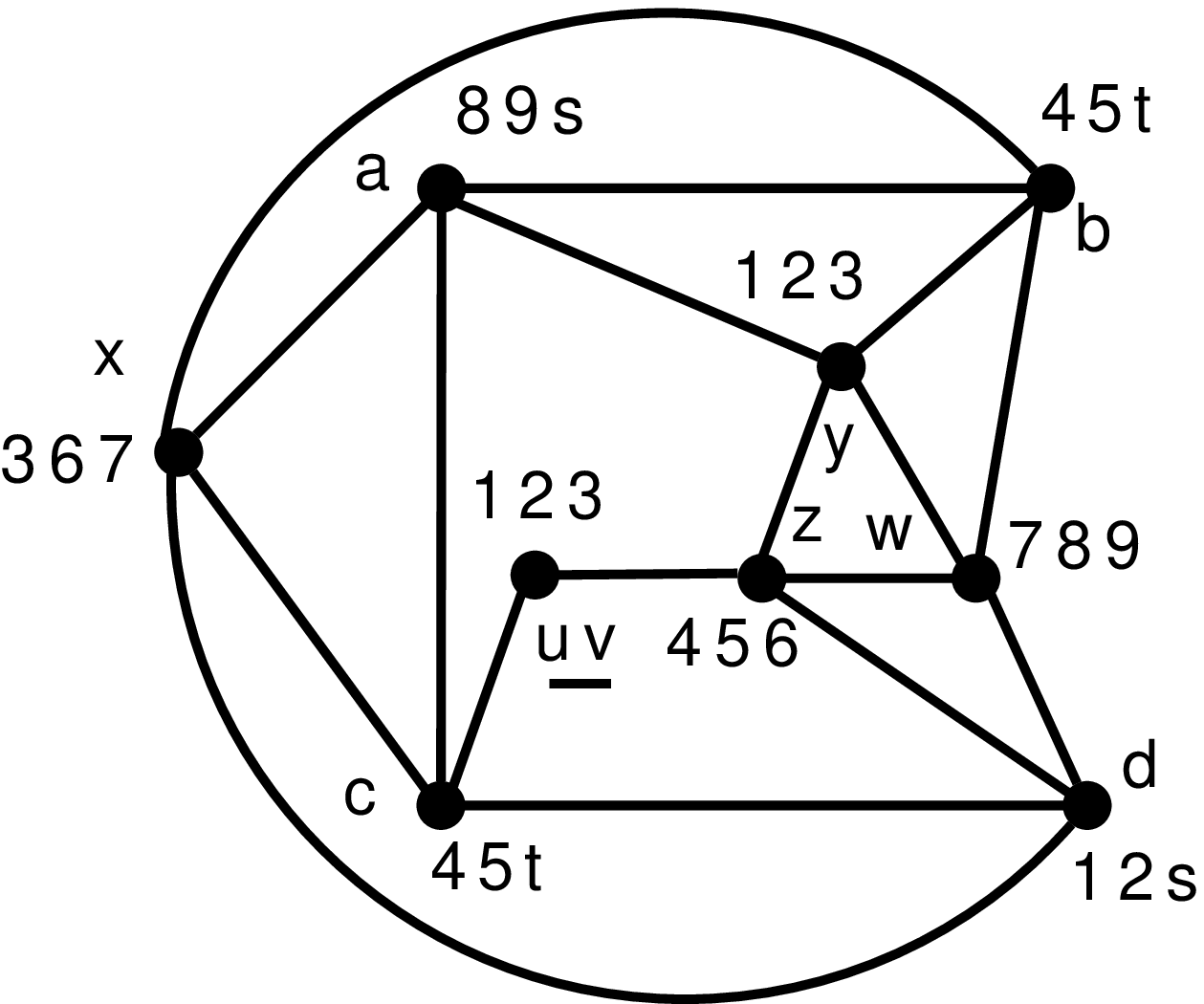, width=0.24\textwidth}
\centerline{$G_2+uv$  \hspace*{1cm} \hfil $G_2/uv$  }
  \caption{Case $H_4$: both graph $G_2+uv$ and $G_2/uv$ are 11:3-colorable.}
  \label{fig:H4}
\end{figure}

Note that $H_4$, $H_6$, and $H_8$ are isomorphic to each other.
It suffices to show $G$ does not contain $H_4$.
Suppose that $H_4$ is a proper induced subgraph of $G$.
Let $G_1$ be the induced subgraph of $G$ by deleting
all vertices in $H_4$. Note $C_8^2$ is not a proper subgraph of any graph in $\G_4$.
We have $G_1\not=C_8^2$.
Note that $c$ and $z$ have degree $3$ while other vertices in $H_4$ have degree $4$.
 Since $G$ is 2-connected, $c$ has a unique neighbor, denoted by $u$,
in $V(G_1)$. Similarly, $z$ has a unique neighbor, denoted by $v$,
in $V(G_1)$. Observe that the pair $\{u,v\}$ forms a vertex cut of $G$.
Let $G_2$ be the induced graph of $G$ on $V(H_4)\cup \{u,v\}$.
Applying Lemma \ref{cut2} to $G$ with $G_1$ and $G_2$, we have
$$\chi_f(G) \leq \max\{\chi_f(G_1), \chi_f(G_2+uv), \chi_f(G_2/uv)\}.$$
Figure \ref{fig:H4} shows  $\chi_f(G_2+uv)$ and $\chi_f(G_2/uv)$ are at most $11/3$.

Since $\chi_f(G)> 11/3$, we have $\chi_f(G)\leq \chi_f(G_1)$.
Now $G_1$ is $K_4$-free and has maximum degree at most $4$; $G_1$ has
fewer vertices than $G$. This contradicts the minimality of $G$.

If $G=H_4$, then $\chi_f(H_4)\leq 11/3$, since $H_4$ is a subgraph of $G_2+uv$
 in Figure \ref{fig:H4}.

\begin{figure}[htbp]
 \centerline{ \psfig{figure=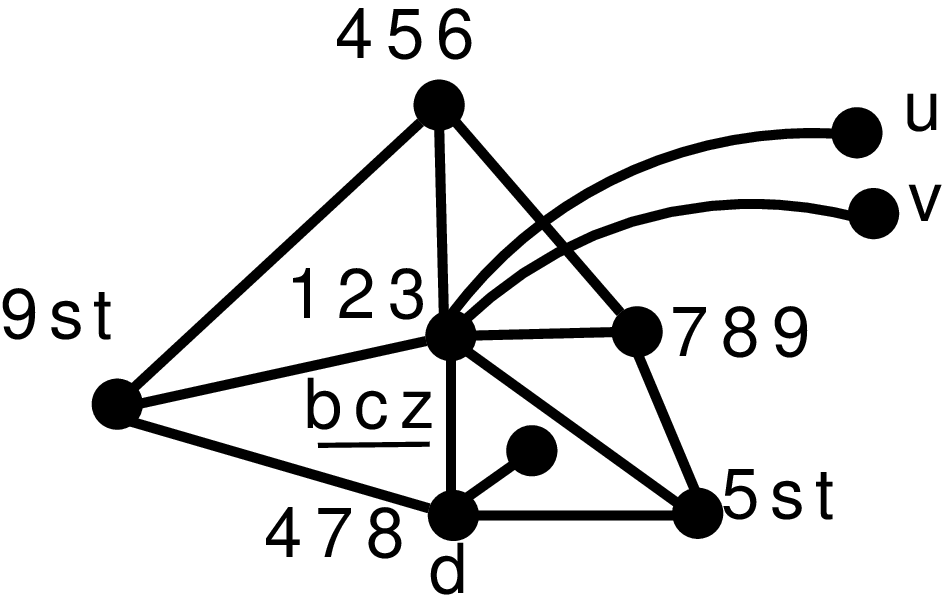, width=0.33\textwidth} \hfill \psfig{figure=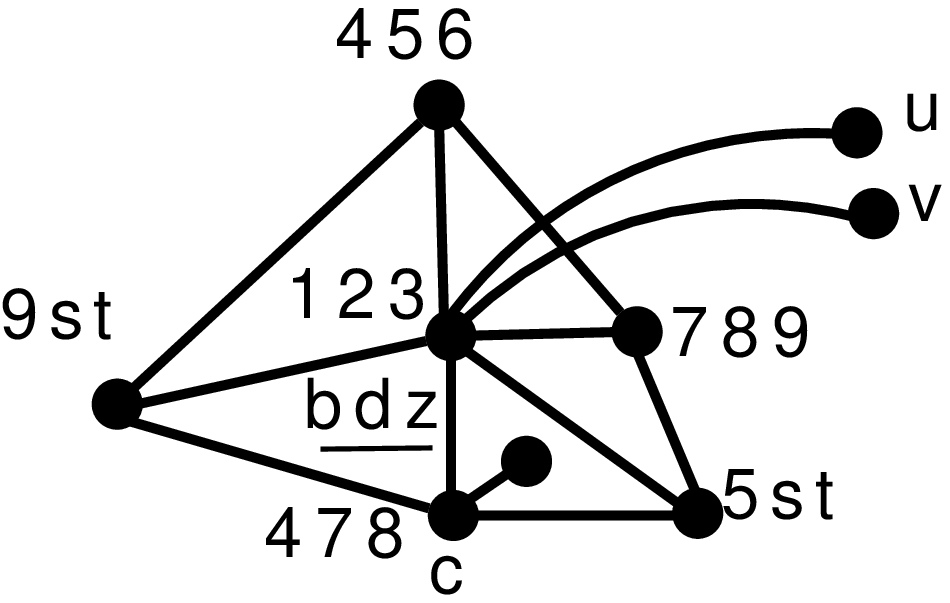,
 width=0.33\textwidth}  \hfill \psfig{figure=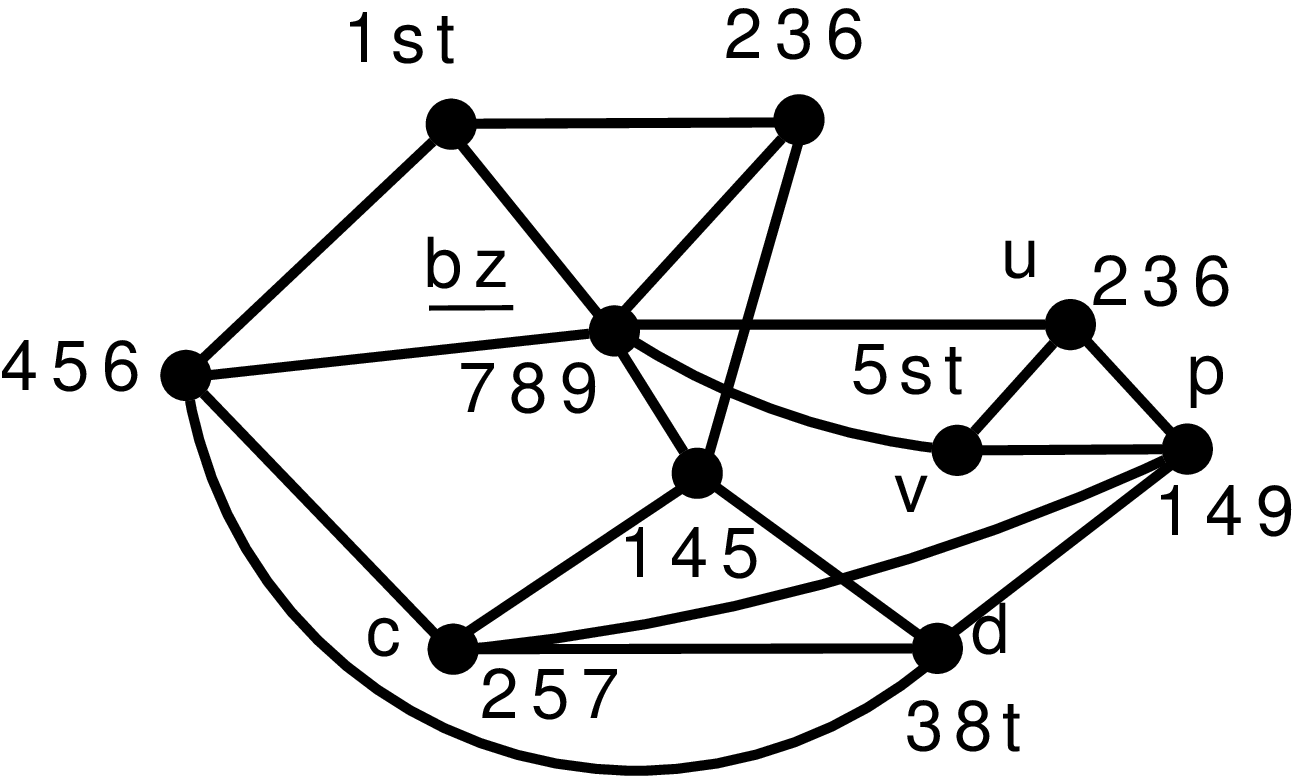,width=0.33\textwidth}}
\centerline{$G_4$ \hspace{3.4cm} $G_5$ \hspace{3.4cm} $G_6$}
\caption{Case $H_9$: the graphs $G_4$, $G_5$, and $G_6$ are 11:3-colorable.}
\label{fig:H9}
\end{figure}

Now we consider the last case $H_9$. First, we contract $b,c,z$ into a
fat vertex denoted by $\underline{bcz}$.  We write $G/bcz$ for the graph after
this contraction. Observe that $\{\underline{bcz},d\}$ is a vertex-cut of
$G/bcz$. Let $G_4$ and $G_4'$ be two connected subgraphs of $G/bcz$
such that $G_4 \cup G_4'=G/bcz$, $G_4 \cap G_4'=\{\underline{bcz},d\}$, and
$\{u,v\} \subset G_4'$. Note that $G_4$ is $11$:$3$ colorable, see
Figure \ref{fig:H9}.  Now by Lemma \ref{cut2}, we have
$$
\chi_f(G/bcz) \leq \max \{\chi_f(G_4), \chi_f(G_4')\}.
$$
As $\{b,c,z\}$ is an independent set,  each $a$:$b$-coloring of
$G/bcz$ gives an $a$:$b$-coloring of $G$, that is $\chi_f(G/bcz)\geq
\chi_f(G)>11/3$. The graph $G_4$ is 11:3-colorable; see Figure
\ref{fig:H9}. Thus we have $\chi_f(G_4')\geq \chi_f(G/bcz)\geq
\chi_f(G)$.  It is easy to check that $G_4'$ has maximum degree $4$,
$K_4$-free,  and it is not $C_8^2$. Hence $G_4'$ must contain a
$K_4$. Otherwise, it contradicts the minimality of
$G$.

Second, we contract $\{b,d,z\}$ into a fat vertex $\underline{bdz}$
and denote the graph by $G/bcz$. Let $G_5$ and $G_5'$ be two
connected subgraphs of $G/bdz$ such that $G_5 \cup G_5'=G/bzd$, $G_5
\cap G_5'=\{\underline{bzd},c\}$, and $\{u,v\} \subset G_5'$. Note
that $G_5$ is 11:3-colorable; see Figure \ref{fig:H9}. By a similar
argument, $G_5'$ must contain a $K_4$.

The remaining case is that both $G_4'$ and $G_5'$ have a $K_4$ when we contract $b$ and $z$.
Since the original graph $G$ is $K_4$-free, the $K_4$ in $G_4'$ (and in
$G_5'$) must contain the fat vertex $\underline{bcz}$ (or
$\underline{bdz}$), respectively. Note that each of the four
vertices $b,c,d,z$ has at most one edge leaving $H_9$. There must be
a triangle $uvp$ in $G$ and these four outward edges are connected
to some element of $\{u,v,p\}$. The graph $G/{bz}$ must contain the
subgraph $G_6$ as drawn in Figure \ref{fig:H9}.

Note that $\{u,v\}$ is a vertex-cut in $G/ bz$. Let $G_6$ and $G_6'$
be two connected subgraphs of $G/bz,$ which satisfy $G_6 \cup
G_6'=G$, $G_6 \cap G_6'=\{u,v\}$, and $\underline{bz} \in G_6$. By
Lemma \ref{cut2}, we have
$$
\chi_f(G/bz) \leq \max\{\chi_f(G_6), \chi_f(G_6')\}.
$$
Note that $G_6$ is 11:3-colorable; see Figure \ref{fig:H9}.  We also
have $\chi_f(G/bz) \geq \chi_f(G) >\frac{11}{3}$. We obtain
$\chi_f(G_6') \geq \chi_f(G/bz) \geq \chi_f(G)$. Observe that $G_6'$
is a subgraph of $G$.  We arrive at a contradiction of the minimality
of $G$.

If $\Gamma_G(x)=C_4$, then the only possible choice for the two
independent sets are $\{a,c\}$ and $\{b,d\}$. If there is some
triangle $H$ such that $(\{a,c\},\{b,d\},H)$ is a bad triple, then
we have
$$
|E(\Gamma_G(x),H)| \geq 5.
$$
However, $|E(\Gamma_G(x),H)| \leq 4$. This is a contradiction. Thus
the lemma follows in this case.

 We can select two
vertex disjoint non-edges $S_1$ and $S_2$ such that the graph
$G/S_1/S_2$ contains no $K_5^-$. For these particular $S_1$ and
$S_2$, if $G/S_1/S_2$ contains no $G_0$, then Lemma \ref{D41} holds.

Without loss of generality, we assume that $G/S_1/S_2$ does contain
$G_0$. Let $s_i=\underline{S_i}$ for $i=1,2$.
 Observe that both $s_1$ and $s_2$ have
four neighbors $u,v,p,q$ other than $x$ in $G_0$.  It follows
that  $$|E(S_1\cup S_2, \{u,v,p,q\})|\geq 8.$$
On the one hand, we have
\begin{eqnarray*}
  |E(G\mid_{S_1\cup S_2})| &=&\frac{1}{2}\left(\sum_{v\in S_1\cup S_2}d(v)
-|E(S_1\cup S_2, \{u,v,p,q\})| - 4\right)\\
&\leq& \frac{1}{2}(16- 8 - 4)\\
&=&2.
\end{eqnarray*}
On the other hand, $\alpha(\Gamma(x))=2$ implies $G\mid_{S_1\cup
S_2}$ contains at least two edges. Thus, we  have $\Gamma_G(x)=2 ~
e$. Label the vertices in $\Gamma_G(x)$ by $a,b,c,d$ as in Figure
\ref{Delta=4}. We assume $ab$ and $cd$ are edges while $ac, bd, ad,
bc$ are non-edges. Observe that each vertex in $\{u,v,p,q\}$ has
exactly two neighbors in $\{a,b,c,d\}$.

If one vertex, say $u$, has two neighbors forming a non-edge, say $ac$, then
we can choose $S'_1=\{a,c\}$ and $S_2'=\{b,d\}$. It is easy to check that
  $G/S_1'/S_2'$ contains neither $G_0$ nor $K_5^-$. We are done in this case.

\begin{figure}[htbp]
\centerline{\psfig{figure=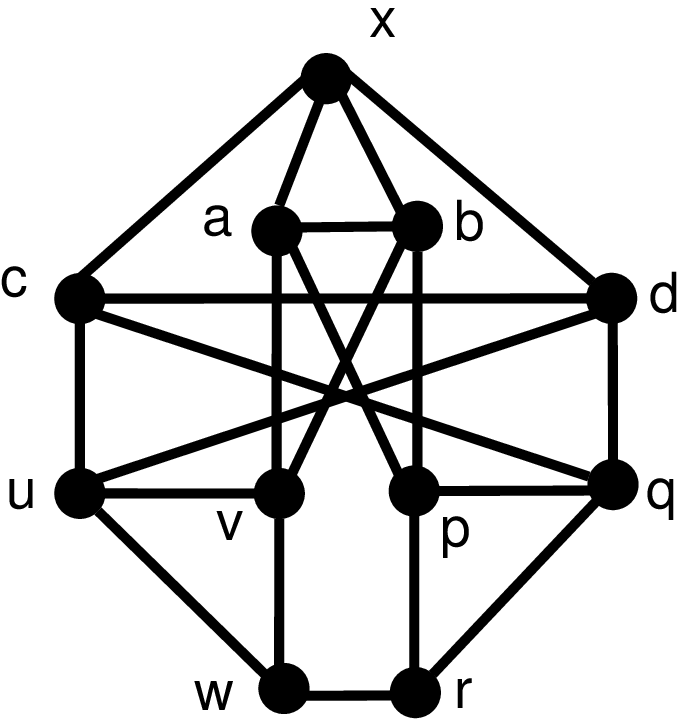, width=0.28\textwidth} \hfil
\psfig{figure=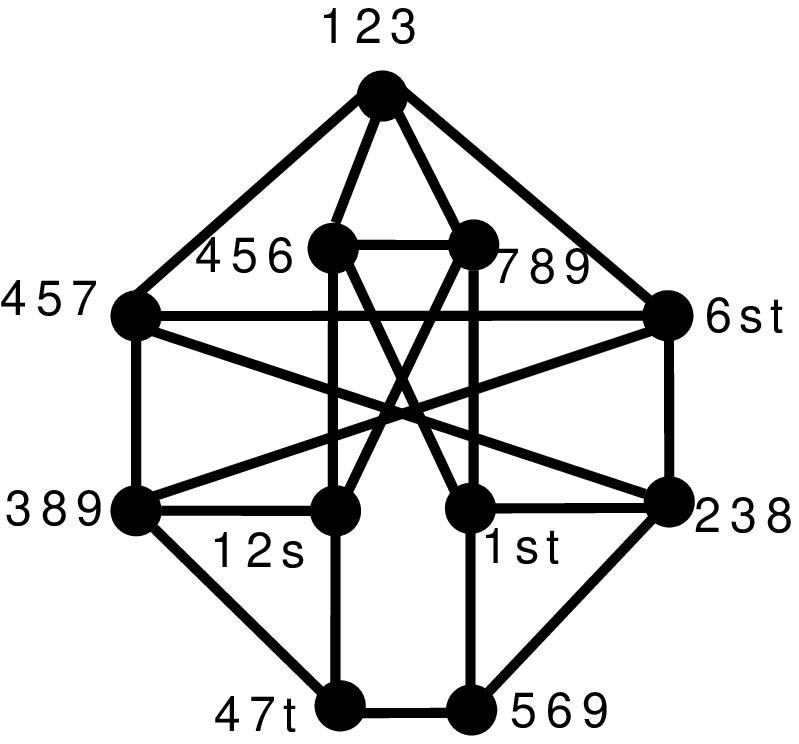, width=0.32\textwidth}}
\caption{$H_{10}$ and an 11:3-coloring of $H_{10}$.} \label{fig:H14}
\end{figure}

In the remaining case, we can assume that for each vertex $y$ in
$\{u,v,p,q\}$, the neighbors of $y$ in $\{a,b,c,d\}$ always form an
edge. Up to relabeling vertices, there is only one arrangement for
edges between $\{u,v,p,q\}$ and $\{a,b,c,d\}$; see the graph
$H_{10}$ defined in Figure \ref{fig:H14}.
The graph $H_{10}$ is 11:3-colorable as shown in Figure
\ref{fig:H14}. Since $\chi_f(G)>11/3$, $H_{10}$ is a proper subgraph
of $G$.
 Note in $H_{10}$, every vertices except $w$ and $r$
has degree $4$; both $w$ and $r$ have degree $3$. Thus, $\{w,r\}$ is
a vertex cut of $G$. Let $G_1=H_{10}$ and $G_2$ be the subgraph of
$G$ by deleting vertices in $\{x,a,b,c,d,p,q,u,v\}$. Applying Lemma
\ref{cut2} with $G_1$ and $G_2$ defined above, we have
$$
\chi_f(G) \leq \max \{\chi_f(G_1),\chi_f(G_2)\}.
$$
Since $\chi_f(G)> 11/3$  and $\chi_f(G_1) \leq 11/3 $ (see Figure
\ref{fig:H14}), we must have $\chi_f(G_2) \geq \chi_f(G)$.  Note
that $G_2$ has fewer number of vertices than $G$. This contradicts
the minimality of $G$. Therefore, the lemma follows. \hfill
$\square$


\end{document}